\documentclass[11pt]{article}
\usepackage{amssymb}

\begin{document}

\begin{center}
\textbf{MEAN CURVATURE MOTION OF GRAPHS WITH CONSTANT CONTACT ANGLE AT A FREE BOUNDARY}\\
by A. Freire\\

University of Tennessee, Knoxville \end{center}

\begin{abstract}
We consider the motion by mean curvature of an $n$-dimensional graph
over a time-dependent domain in $\mathbb{R}^n$, intersecting
$\mathbb{R}^n$ at a constant angle. In the general case, we prove
local existence for the corresponding quasilinear parabolic equation
with a free boundary, and derive a continuation criterion based on
the second fundamental form. If the initial graph is concave, we
show this is preserved, and that the solution exists only for finite
time. This corresponds to a symmetric version of mean curvature
motion of a network of hypersurfaces with triple junctions, with
constant contact angle at the junctions.
\end{abstract}


\vspace{.5cm}

\textbf{1. Time-dependent graphs with a contact angle condition.}
\vspace{.2cm}

We consider a moving hypersurface $\Sigma_t$ in $\mathbb{R}^{n+1}$,
with normal velocity equal to its mean curvature, assumed to be a
graph over a time-dependent open set $D(t)\subset {\mathbb R}^n$
(not necessarily bounded, or connected.) The (properly embedded)
$(n-1)$-submanifold of intersection:
$$\Gamma(t)=\Sigma_t\cap {\mathbb R^n}=\partial D(t)$$
is a `moving boundary'. Along $\Gamma(t)$ we impose a constant-angle
condition:
$$\langle N,e_{n+1}\rangle_{|\Gamma(t)}=\beta,$$
where $0<\beta<1$ is a \emph{constant} and $N$ is the upward unit
normal of $\Sigma_t$. `Mean curvature motion' is defined by the law:
$$V_N=H,$$
where $V_N=\langle V,N\rangle$, with $V=\partial_tF$ the velocity
vector in a given parametrization $F(t)$ of $\Sigma_t$ ($V$ depends
on the parametrization, while $V_N$ does not). A particular
parametrization yields `mean curvature flow':
$$\partial_tF=HN.$$
For graphs, it is natural to consider `graph mean curvature motion':
if $\Sigma_t=\mbox {graph }w(t)$ for a function
$w(t):D(t)\rightarrow {\mathbb R}$, imposing $\langle
\partial_tF,N\rangle=H$ with $F(y,t)=[y,w(y,t)]$ for $y\in D(t)$, we
find:
$$w_t=\sqrt{1+|D w|^2}H$$
(and the velocity is vertical, $\partial_tF=w_te_{n+1})$. With the
contact angle condition, we obtain a free boundary problem for a
quasilinear PDE:

$$\left\{ \begin{array}{l}w_t=g^{ij}(Dw)w_{ij}\quad \mbox{ in } D(t),\\
w=0,\quad \beta\sqrt{1+|Dw|^2}=1\mbox{ on }\partial
D(t),\end{array}\right .$$ where $g^{ij}
(Dw)=\delta^{ij}-w_iw_j/(1+|Dw|^2)$ is the inverse metric
matrix.\vspace{.2cm}

\emph{Remark 1.1.} It is easy to see that the constant-angle
boundary condition is \emph{incompatible} with mean curvature flow
parametrized over a \emph{fixed} domain $D_0$: on $\partial D_0$ we
would have $\langle F,e_{n+1}\rangle=0$, leading to $\langle
\partial_tF,e_{n+1}\rangle=0$, incompatible with $\partial_tF=HN$ and
$\langle N,e_{n+1}\rangle=\beta$. If we parametrize over a
time-dependent domain, mean curvature flow leads to a normal
velocity for the moving boundary that is difficult to control; hence
we chose to analyze the geometry of the motion in terms of the
`graph m.c.m.' parametrization.\vspace{.2cm}

 To establish short-time existence (in parabolic H\"{o}lder
spaces) we will work with a third parametrization of the motion,
defined over a fixed domain:
$$F(t):D_0\rightarrow {\mathbb R}^{n+1},\quad
F(x,t)=[\varphi(x,t),u(x,t)]\in {\mathbb R}^n\times {\mathbb R},$$
where $\varphi(t):D_0\rightarrow D(t)$ is a diffeomorphism and $F$
is a solution of the parabolic system:
$$F_t=g^{ij}(DF)F_{ij},$$
where $g_{ij}=\langle F_i,F_j\rangle$ is the induced metric on
$\Sigma_t$ and $g^{ij}$ is the inverse metric matrix. \vspace{.2cm}

In the first part of the paper (sections 3 to 8) we prove the
following short-time existence theorem (on $Q:=D_0\times
[0,T]$):\vspace{.2cm}

\textbf{Theorem 1.1.} Let $\Sigma_0\subset {\mathbb R}^{n+1}$ be a
$C^{3+\bar{\alpha}}$ graph over $D_0\subset {\mathbb R}^n$
satisfying the contact and angle conditions at $\partial D_0$. There
exists a parametrization $F_0=[\varphi_0,u_0]\in C^{2+\alpha}(D_0)$
of $\Sigma_0$ (where $\alpha=\bar{\alpha}^2$ and $\varphi_0$ is a
diffeomorphism of $D_0$ satisfying the `orthogonality conditions' on
$\partial D_0$), $T>0$ depending only on $\Sigma_0$ and a unique
solution $F\in C^{2+\alpha,1+\alpha/2}(Q^T;{\mathbb R}^{n+1})$ of
the system:
$$\left \{ \begin{array}{l}
\partial_t F-g^{ij}(DF)\partial_i\partial_jF=0,\quad F=[\varphi,u]\in \mathbb{R}^n\times \mathbb{R}\\
u_{|\partial D_0}=0, \quad N^{n+1}(DF)_{|\partial
D_0}=\beta,\end{array}\right .$$ with initial data $F_0$ and where
$\varphi(t):D_0\rightarrow D(t)\subset \mathbb{R}^n$ is a
diffeomorphism satisfying, in addition, the `orthogonality
conditions' at $\partial D_0$ (described in section 3.)
\vspace{.2cm}

The system and boundary conditions are discussed in more detail in
section 3. Sections 4, 5, and 6 deal with compatibility at $t=0$,
linearization and the verification that the boundary conditions
satisfy `complementarity'. In particular, adjusting the initial
diffeomorphism $\varphi_0$ to ensure compatibility (section 4) leads
to the `loss of differentiability' seen in theorem 1.1. The required
estimates in H\"{o}lder spaces for the linearized system are
described in section 7, and the proof concluded (by a fixed-point
argument) in section 8. While the general scheme is standard, since
we are dealing with a free boundary problem with somewhat
non-standard boundary conditions, details are included. Free
boundary-type problems for mean curvature motion of graphs have
apparently not previously been considered.

We describe the evolution equations in the rotationally symmetric
case in section 9 (including a stationary example for the exterior
problem) and the extension to the case of a graph motion $\Sigma_t$
intersecting fixed support hypersurfaces orthogonally in section 10.
\vspace{.2cm}

The original motivation for this work was to establish (by classical
parabolic PDE methods) existence-uniqueness for mean curvature
motion of networks of surfaces meeting along triple junctions with
constant-angle conditions. One can use a motion $\Sigma_t$ of graphs
with constant contact angle to produce examples of `triple junction
motion': three hypersurfaces moving by mean curvature meeting along
an ($n-1$)-dimensional submanifold $\Sigma(t)$ so that the three
normals make constant angles (say, 120 degrees) along $\Gamma(t)$.
We simply reflect on $\mathbb R^n$, so the hypersurfaces are
$\Sigma_t, \bar{\Sigma_t}$ and $\mathbb{R}^n-\bar{D}(t)$. If
$\Sigma_t=\mbox{graph }w(t)$ with $w>0$, the system is
\emph{embedded} in ${\mathbb R}^{n+1}$. This is mean curvature
motion of a `symmetric triple junction of graphs'.\vspace{.2cm}

Short-time existence holds for general triple junctions of graphs
moving by mean curvature with constant 120 degree angles at the
junction, provided a compatibility condition holds along the
junction (see section 16). Since the free-boundary problem is easier
to understand in the symmetric case, we decided to do this first. In
addition, in the present case it is possible to go further towards a
geometric global existence result. Motivated by recent work on
`lens-type' curve networks \cite{LensSeminar}, in the second part of
the paper (sections 11-15) we consider continuation criteria and
preservation of concavity. Since we chose to develop these results
for graph motion with a free boundary, although the general lines of
proof (via maximum principles) have precedents, the details of the
arguments are new. For example, section 12 contains an extension of
the maximum principle for symmetric tensors with Neumann-type
boundary conditions given in \cite{Stahl}, which in our setting
allows one to show preservation of weak concavity in general.
Section 14 includes a continuation criterion for the flow. The
results obtained in sections 11-15 are summarized in the following
theorem. ($h$ denotes the second fundamental form, pulled back to a
symmetric 2-tensor on $D(t)$.)\vspace{.2cm}

\textbf{Theorem 1.2.} If $\Sigma_0$ is weakly concave ($h\leq 0$ at
$t=0$), this is preserved by the evolution. Let $T_{max}$ be the
maximal existence time for the evolution. If the mean curvature of
$\Sigma_0$ is strictly negative ($\sup_{\Sigma_0}H=H_0<0$), then
$T_{max}$ is finite. Assuming $T_{max}<\infty$, we have:
$$\limsup_{t\rightarrow T_{max}}[\sup_{\Gamma_t}(|h|_g+|\nabla^{tan}h^{tan}|_g)]=\infty$$
(if $n=2$, in the concave case). If there is no `gradient blowup' at
$T_{max}$, the hypersurface contracts to a compact convex subset of
$\mathbb{R}^n$ as $t\rightarrow T_{max}$.\vspace{.2cm}

\emph{Remark 1.2.} We have not yet proved that the diameter tends to
zero as $t\rightarrow T_{max}$, though this seems likely based on
the experience with curves \cite{LensSeminar} (in the absence of
`gradient blowup'.) It is an interesting question (even in the
concave case, for $n=2$) whether the latter can really occur (that
is, $\sup_{\Gamma_t}|\nabla^{tan}h^{tan}|_g\rightarrow \infty$ as
$t\rightarrow T_{max}$, while $|h|_g$ remains bounded on
$\Gamma_t$.)\vspace{.2cm}

\emph{Acknowledgments.} It is a pleasure to thank Nicholas Alikakos
for originally proposing to consider the problem of mean curvature
flow for networks of surfaces meeting at constant angles, and for
his interest in this work. Most of the work on short-time existence
was done during a stay at the Max-Planck Institute for Gravitational
Physics in Golm (January-June, 2007); I am grateful to the
Max-Planck Society for supporting the visit and to the Geometric
Analysis group at the AEI for the invitation. Finally, thanks to
Mariel S\'{a}ez for communicating the results of the Lens Seminar
(\cite{LensSeminar}) and of her recent work on mean curvature flow
of networks (partly in collaboration with Rafe Mazzeo,
\cite{MazzeoSaez}).\vspace{.4cm}

\textbf{2. Normal velocity of the moving boundary.} The evolution is
naturally supplied with initial data $\Sigma_0$, a graph meeting
$\mathbb{R}^{n+1}$ at the prescribed angle. Since we are interested
in classical solutions in the parabolic H\"{o}lder space
$C^{2+\alpha,1+\alpha/2}$, we expect an additional compatibility
condition at $t=0$. We discuss this first for graph m.c.m. $w(y,t)$.

Denote by $\Gamma(t)$ a global parametrization of $\partial D(t)$
(with domain in a fixed manifold, and `space variables' left
implicit). Differentiating in $t$ the `contact condition'
$w(\Gamma(t),t)=0$, we find:
$$w_t+\langle Dw,\dot{\Gamma}(t)\rangle =0.$$
Denote by $n=n_t$ the unit normal vector field to $\Gamma(t)$,
chosen so that the directional derivative $d_nw>0$. The contact
condition also implies the gradient of $w$ is purely normal:
$$Dw_{|\partial D(t)}=(d_{n}w)n.$$
Combining this with the angle condition, and bearing in mind that
$d_{n}w_{|\Gamma(t)}>0$, we find:
$$d_{n}w=\frac {\beta_0}{\beta}\mbox{ on }\partial D(t),\quad \beta_0:=\sqrt{1-\beta^2}.$$
(In fact, this is a more convenient form of the `angle' boundary
condition for $w$, since it is linear.) Thus, on $\partial D(t)$:
$$\frac 1{\beta}H=\sqrt{1+(d_{n}w)^2}H=w_t=-\langle
\dot{\Gamma}(t),n\rangle d_{n}w=-\dot{\Gamma}_n(t)\frac
{\beta_0}{\beta},$$ and we find the normal velocity of the moving
boundary (independent of the parametrization of $\Gamma_t$):
$$\dot{\Gamma}_n=-\frac 1{\beta_0}H_{|\Gamma(t)},$$
which in particular must hold at $t=0$. Note that we don't get a
`compatibility condition' in the usual sense (of a constraint on the
2-jet of the initial data), but instead an equation of motion for
the moving boundary. (Later, in the fixed-domain formulation, we
will have to deal with a real compatibility condition).
\vspace{.2cm}

\emph{Remark 2.1.} We remark that for more general (non-symmetric,
non-flat) triple junctions with 120 degree angles, the condition:
$$H^1+H^2=H^3\mbox{ on }\Gamma(t)$$
must hold at the junction (for graphs, oriented by the upward
normal), which in particular gives a geometric constraint on the
initial data, for classical evolution in $C^{2+\alpha, 1+\alpha/2}$.
This is automatic in the symmetric case ($w^2=-w^1,w^3\equiv 0$),
since $H^3=0$ and $H^I=tr_{g^I}d^2w^I$ for $I=1,2$. \vspace{.3cm}

\textbf{3. Choice of `gauge'.} It is traditional in moving boundary
problems to parametrize the time-dependent domain $D(t)$ of the
unknown $w(y,t)$ by a time-dependent diffeomorphism:
$$y=\varphi (x,t),\quad \varphi(t):D_0\rightarrow D(t),$$
and then derive the equation satisfied by the coordinate-changed
function from the equation for $w$ (see e.g. \cite{LunardiBaconneau}
or \cite{Solonnikov}). Motivated by the work on curve networks
(\cite{Mantegazza et al.}) we will, instead, consider a general
parametrization:
$$F:D_0\times [0,T]\rightarrow {\mathbb R}^{n+1},\quad
F(x,t)=[\varphi(x,t),u(x,t)]\in\mathbb{R}^n\times \mathbb{R}$$ and
derive an equation for $F$ directly from the definition of mean
curvature motion:
$$\langle \partial_tF,N\rangle=H.$$
(We'll still assume $\varphi(t):D_0\rightarrow D(t)$ is a
diffeomorphism.) The first and second fundamental forms are given
by:
$$g_{ij}=\langle F_i,F_j\rangle, \quad A(F_i,F_j)=\langle
F_{ij},N\rangle.$$ (Notation: $DF=F_ie_i, D^2F(e_i,e_j)=F_{ij}$,
$(e_i)$ is the standard basis of $\mathbb R^{n+1}$.) The mean
curvature is the trace of $A$ in the induced metric:
$$H=\langle g^{ij}(DF)F_{ij},N\rangle.$$ The equation for $F$ is:
$$\langle \partial_tF-g^{ij}(DF)F_{ij},N\rangle=0.$$

There is a natural `gauge choice' yielding a quasilinear parabolic
system:
$$\partial_tF-g^{ij}(DF)F_{ij}=0.$$
We will sometimes refer to this as the `split gauge', since in terms
of the components $F=[\varphi,u]$ we have the essentially decoupled
system:
$$\left \{ \begin{array}{ccc}
\partial_tu-g^{ij}(D\varphi,Du)u_{ij}&=&0,\\
\partial_t\varphi-g^{ij}(D\varphi,Du)\varphi_{ij}&=&0.\end{array} \right .$$ The
splitting is useful to state the boundary conditions:
$$\left \{ \begin{array}{l}u_{|\partial D_0}=0\quad \mbox{ (`contact')},\\
N^{n+1}(D\varphi,Du)_{|\partial D_0}=\beta\quad \mbox{ (`angle')}.
\end{array}\right.$$
 We immediately see there is a problem,
since we have 2 scalar boundary conditions for $n+1$ unknowns (and
no moving boundary to help!) Our solution to this is to introduce
$n-1$ additional `orthogonality conditions' at the boundary for the
parametrization $\varphi(t)$. We impose:
$$\langle D_{\tau}\varphi,D_n\varphi\rangle_{|\partial D_0}=0\quad \mbox{(`orthogonality')},$$
for any $\tau\in T\partial D_0$, where $n$ denotes the inward unit
normal to $D_0$. (We fix a tubular neighborhood $\cal N$ of
$\partial D_0$ and extend $n$ to $\cal N$ so that $d_nn=0$ in $\cal
N$.)

Geometrically, the `orthogonality' boundary condition has precedent
in a method often adopted when dealing with the evolution of
hypersurfaces in $\mathbb{R}^{n+1}$ intersecting a fixed
$n$-dimensional `support surface' orthogonally (see e.g.
\cite{Struwe}): one replaces vanishing inner product of the unit
normals (a single scalar condition) by a stronger Neumann-type
condition for the parametrization, corresponding to $n-1$ scalar
conditions. (More details are given in Section 10.)\vspace{.2cm}

The system must also be supplied with initial data. We assume given
an initial hypersurface $\Sigma_0$, the graph of a
$C^{3+\bar{\alpha}}$ function $\tilde u_0(x)$ defined in the
$C^{3+\bar{\alpha}}$ domain $D_0\subset \mathbb{R}^n$. (The reason
for this choice of differentiability class will be seen later.) It
would seem natural to set $\varphi_0=Id_{D_0}$, but this causes
problems (related to compatibility; see Section 4 below). We do
require the 1-jet of $\varphi_0$ at the boundary to be that of the
identity:
$${\varphi_0}_{|\partial D_0}=Id,\quad D{\varphi_0}_{|\partial D}
=\mathbb{I}.$$ (In particular, the orthogonality condition holds at
$t=0$.)

We need a more explicit expression for the unit normal, and for that
we use the `vector product':
$$\tilde{N}(D\varphi,Du):=(-1)^n\det\left
[\begin{array}{ccc}e_1&\cdots&e_{n+1}\\DF^1&\cdots&DF^{n+1}\end{array}
\right]=(-1)^{n}\det
\left[\begin{array}{cccc}e_1&\ldots&e_n&e_{n+1}\\D\varphi^1&\ldots&D\varphi^n&Du\end{array}\right]$$
$$:=[J(D\varphi,Du),J_{\varphi}]\in {\mathbb R}^n\times
{\mathbb R},$$
 where $DF^i\in {\mathbb R}^n$ for $i=1,\ldots n+1$, $J_{\varphi}>0$ is the jacobian of $\varphi$ and
$(-1)^n$ is introduced to make sure the last component is positive.
$J(D\varphi,Du)$ is an ${\mathbb R}^n$-valued multilinear form,
linear in the components $u_i$ of $Du$ and of weight $n-1$ in the
components of $D\varphi$. It is easy to check that
$J(\mathbb{I},Du)=-Du$. The unit normal is:
$$N(D\varphi,
Du)=\tilde{N}(D\varphi,Du)/(|J(D\varphi,Du)|^2+(J_{\varphi})^2)^{1/2}.$$
Thus the angle condition may be stated in the form:
$$\beta[|J(Du,D\varphi)|^2+(J_{\varphi})^2]^{1/2}_{|\partial
D_0}={J_{\varphi}}_{|\partial D_0},$$ and we lose nothing by
squaring it:
$$B(D\varphi,Du):=
\beta^2|J(Du,D\varphi)|^2-\beta_0^2(J_{\varphi})^2_{|\partial
D_0}=0.$$ \vspace{.3cm}

\textbf{4. Compatibility and the choice of $\varphi_0$.} Assume
$D{\varphi_0}_{|\partial D_0}=\mathbb{I}$. Differentiating in $t$
the contact condition $u_{|\partial D_0}=0$ and evaluating at $t=0$,
we find:
$$0=g^{ij}(\mathbb I,Du_0)u_{0ij}\equiv g_0^{ij}u_{0ij}\mbox{ on
}\partial D_0.$$ To interpret this condition, consider the mean
curvature at $t=0$, on $\partial D_0$:
$$H_0=\frac 1{v_0}[\langle
J(\mathbb{I},Du_0),g_0^{ij}\varphi_{0ij}\rangle+J_{\varphi_0}g_0^{ij}u_{0ij}],$$
where:
$$v_0=[|J(\mathbb{I},Du_0)|^2+J_{\varphi_0}^2]^{1/2}_{|\partial D_0}
=(|Du_0|^2+1)^{1/2}_{|\partial D_0}=\frac 1{\beta},$$ using (recall
$\beta_0:=\sqrt{1-\beta^2}$):
$$J(\mathbb I,Du_0)=-Du_0=-(D_nu_0)n=-\frac
{\beta_0}{\beta}n$$ on $\partial D_0$. Thus the compatibility
condition is equivalent to:
$${H_0}_{|\partial D_0}=-\beta_0 g_0^{ij}\langle
\varphi_{0ij},n\rangle_{|\partial D_0}.$$ This implies we can't
choose $\varphi_0\equiv Id$ (on all of $D_0$), unless
${H_0}_{|\partial D_0}\equiv 0$, a constraint not present in the
geometric problem (as seen above).\footnote{The compatibility
condition ${H_0}_{|\partial D_0}=0$ does occur for graph m.c.m. with
Dirichlet boundary conditions in a mean-convex domain
\cite{Huisken89}.} Instead, regarding $H_0$ as given (by
$\Sigma_0$), and using:
$$g_0^{ij}=\delta_{ij}-\frac{u_{0i}u_{0j}}{v_0^2}=\delta_{ij}-\beta_0^2n^in^j,$$
we find the compatibility constraint:
$$\langle
(\delta_{ij}-\beta_0^2n^in^j)\varphi_{0ij},n\rangle=-\frac
1{\beta_0}H_0\mbox{ on }\partial D_0.$$ Given the zero and first
order constraints on $\varphi_0$, this can also be written as:
$$n^in^j\langle \varphi_{0ij},n\rangle=-\frac
1{\beta^2\beta_0}H_0\mbox{ on }\partial D_0.$$ The next lemma shows
this can be solved.\vspace{.2cm}

\textbf{Lemma 4.1.} Let $D_0\subset \mathbb{R}^n$ be a uniformly
$C^{3+\alpha}$ domain (possibly unbounded), $h\in
C^{\alpha}(\partial D_0)$ $(0<\alpha<1)$.

(i) One may find a diffeomorphism $\varphi\in Diff^{2+\alpha}(D_0)$
satisfying on $\partial D_0$:
$$\varphi=Id, \quad d\varphi=\mathbb{I},\quad n\cdot
d^2\varphi(n,n)=h.$$

(ii) More generally, given a non-vanishing vector field $e\in
C^{1+\alpha}(\partial D_0;{\mathbb R}^n)$, one may find $\varphi\in
Diff^{2+\alpha}(D_0)$ satisfying on $\partial D_0$:
$$\varphi=Id, \quad d_n\varphi=e, \quad n\cdot
d^2\varphi(n,n)=h.$$ If $\partial D_0$ has two components, we may
even require $\varphi$ to satisfy the conditions in parts (i) and
(ii) at $\partial_1D_0$, $\partial_2D_0$ (resp.), with different
functions $h$. (This will be needed in section 10).
 \vspace{.1cm}

As usual, a domain is `uniformly $C^{3+\alpha}$' if at each boundary
point there are local charts to the upper half-space (of class
$C^{3+\alpha}$), defined on balls of uniform radius, and with
uniform bounds on the $C^{3+\alpha}$ norms of the charts and their
inverses.\vspace{.2cm}

\emph{Remarks:} 4.1.The proof is given in Appendix 1.

4.2. Note that, in particular, $\varphi$ satisfies the orthogonality
conditions at $\partial D_0$.

4.3. It is at this step in the proof that we have a drop in
regularity: for $C^{2+\alpha}$ local solutions, we require
$C^{3+\alpha}$ initial data. While this is not unexpected in
free-boundary problems (see e.g. \cite{LunardiBaconneau}), I don't
know a counterexample to the lemma if $D_0$ is assumed to be a
$C^{2+\alpha}$ domain.

4.4. In our application of the lemma, we in fact have $h\in
C^{1+\alpha}(\partial D_0)$, but this does not imply higher
regularity for $\varphi$. \vspace{.3cm}

\textbf{5. Linearization.} The evolution equation and boundary
conditions in `split gauge' are: $$\left \{\begin{array}{ccc}
 F_t-g^{ij}(DF)F_{ij}&=&0,\\
 u_{|\partial D_0}&=&0,\\
 B(D\varphi,Du)_{|\partial
D_0}&=&0,\\
{\cal O}(D\varphi)_{|\partial D_0}&=&0,
\end{array}\right .$$
where: $$ {\cal O}(D\varphi):=\langle
D^T\varphi,D_n\varphi\rangle.$$ Here
$D^T\varphi=D\varphi-(d_n\varphi)\langle \cdot,n\rangle$ is an
${\mathbb R}^n$-valued $(n-1)$-form on $\partial D_0$. We'll prove
short-time existence for this system (with initial data
$u_0,\varphi_0$) in $C^{2+\alpha, 1+\alpha/2}$ by the usual
fixed-point argument based on linear parabolic theory. Given
$\bar{F}=[\bar{\varphi},\bar{u}]$ in a suitable ball in this
H\"older space with center $F_0=[\varphi_0,u_0]$, it suffices to
consider the `pseudolinearization' of the system:
$$F_t-g^{ij}(DF_0)F_{ij}=[g^{ij}(D\bar{F})-g^{ij}(DF_0)]\bar{F}_{ij}:={\cal
F}(\bar{F},{F_0}):=\bar{\cal F}\quad (LPDE);$$ a fixed point of the
map $\bar{F}\mapsto F$ corresponds to a solution of the quasilinear
equation.\vspace{.2cm}

For the nonlinear boundary conditions, we need the honest
linearization at $F_0$. For the angle condition, a computation using
the boundary constraints on $u_0$ and $\varphi_0$ yields:
$$\frac 12 {\cal
L}_0B[D\varphi,Du]=-\beta\beta_0d_nu-\beta_0^2\langle
d_n\varphi,n\rangle.$$ The corresponding linear boundary condition
will be:
$$\beta\beta_0d_nu+\beta_0^2\langle
d_n\varphi,n\rangle={\cal B}(D\bar F,DF_0):=\bar{\cal B},$$ where:
$$2{\cal B}(DF^1,DF^2):=B(D\varphi_1,Du_1)-B(D\varphi_2,Du_2)-{\cal
L}_0B[D(\varphi_1-\varphi_2),D(u_1-u_2)],$$ and we used:
$$-\frac 12 {\cal L}_0[D\varphi_0,Du_0]_{|\partial
D_0}=\beta\beta_0d_nu_0+\beta_0^2\langle
d_n\varphi_0,n\rangle_{|\partial D_0}=0.$$ Also,
$B(D\varphi_0,Du_0)_{|\partial D_0}=0$, so at a fixed point
$B(D\varphi,Du)_{|\partial D_0}=0$.\vspace{.2cm}

Linearizing the orthogonality boundary condition, we find that
${\cal L}_0{\cal O}[D\varphi]$ is the $(n-1)$-form on $\partial
D_0$:
$${\cal
L}_0{\cal
O}[D\varphi](v)=(\partial_j\varphi^i+\partial_i\varphi^j)n^j(\delta_{ik}-n^kn^i)v^k$$
(with sum over repeated indices.) The corresponding linear boundary
condition is:
$$\langle d_n\varphi,proj^T\rangle+\langle
D^T\varphi,n\rangle=-\Omega(D\bar{\varphi},D\varphi_0):=\bar{\Omega},$$
where $proj^T$ denotes orthogonal projection
$\mathbb{R}^n\rightarrow T\partial D_0$, and:
$$\Omega(D\varphi_1,D\varphi_2):={\cal O}(D\varphi_1)-{\cal O}(D\varphi_2)-{\cal
L}_0{\cal O}[D\varphi_1-D\varphi_2],$$ and we used:
$${\cal L}_0{\cal O}[D\varphi_0]_{|\partial D_0}=\langle
(d_n\varphi_0)^T,\cdot\rangle+\langle
D^T\varphi_0,n\rangle_{|\partial D_0}=0.$$ \vspace{.3cm}

\textbf{6. Complementarity.} We wish to apply linear existence
theory to the system:
$$F_t-g^{ij}(DF_0)F_{ij}=\bar{\cal F},$$
with boundary conditions at $\partial D_0$:
$$(LBC)\left \{ \begin{array}{l}u=0\\
\beta\beta_0d_nu+\beta_0^2\langle d_n\varphi,n\rangle =\bar{\cal
B},\\
\langle d_n\varphi,proj^T\rangle +\langle
D^T\varphi,n\rangle=-\bar{\Omega}\end{array}\right .$$ and initial
conditions:
$$u_{t=0}=u_0,\quad \varphi_{t=0}=\varphi_0.$$
It is easy to see that the initial data satisfy the linearized
boundary conditions, and above we constructed $\varphi_0$ so as to
guarantee $g^{ij}(Du_0,D\varphi_0)u_{0ij}{_{|\partial D_0}}=0.$
(There is no first-order compatibility condition for $\varphi_0$.)
Thus the linear system satisfies the required compatibility at
$t=0$.\vspace{.2cm}

Since the linearized boundary conditions are slightly non-standard,
we must verify they satisfy the `complementarity'
(Lopatinski-Shapiro) conditions. We \emph{fix} $x_0\in \partial D_0$
and introduce adapted coordinates $(\rho,\sigma)$ in a neighborhood
${\cal N}_0\subset {\cal N}$ of $x_0$ in $D_0$:
$$x\in {\cal N}_0\Rightarrow x=\Gamma_0(\sigma)+\rho n(\sigma),\quad \sigma=(\sigma_a)\in {\cal U},$$
where $\Gamma_0:{\cal U}\rightarrow {\mathbb R^n}$ is a local chart
for $\partial D_0$ at $x_0$ (${\cal U}\subset {\mathbb R}^{n-1}$
open). This defines a basis of tangential vector fields in
$\Gamma_0({\cal U})$, and we may assume that, at $x_0$: $\langle
\tau_a,\tau_b\rangle=\delta_{ab}$ and $\nabla_{\tau_a}\tau_b(x_0)=0$
(for the induced connection on $T\partial D_0$). Let $U$ and $\psi$
be defined in $(-\rho_1,0)\times {\cal U}\times [0,T]$ by:
$$U(\rho,\sigma,t)=u(\Gamma_0(\sigma)+\rho n(\sigma),t),\quad
\psi(\rho,\sigma,t)=\varphi(\Gamma_0(\sigma)+\rho n(\sigma),t).$$ In
these coordinates, the induced metric is written (in `block form'):
$$[g(DF_0)]=\left [ \begin{array}{cc}|\psi_{\rho}|^2+(U_{\rho})^2 &\langle
\psi_{\rho},\psi_{a}\rangle+U_{\rho}U_a\\\langle
\psi_{\rho},\psi_a\rangle+U_{\rho}U_a &\langle
\psi_{a},\psi_b\rangle+U_aU_b\end{array}\right ]_{|t=0}=\left
[\begin{array}{cc}\frac 1{\beta^2}& 0\\0& \mathbb{I}_{n-1}
\end{array}\right ]$$ at $t=0$ and $x_0$.

We have: $$U_{\rho \rho}=D^2u(n,n)\mbox{ (since $d_nn=0$)},$$

$$U_{ab}=D^2u(\tau_a,\tau_b)+Du\cdot
\nabla_{\tau_a}\tau_b=D^2u(\tau_a,\tau_b) \mbox{ at }x_0,$$ and we
don't need $U_{\rho a}$, since $g_{\rho a}=0$ at $x_0$.

Thus: $$tr_{g_0}D^2u(x_0)=\beta^2 D^2u(n,n)+\sum_a
D^2u(\tau_a,\tau_a)=\beta^2U_{\rho \rho}+\sum_aU_{aa}:=
\beta^2U_{\rho \rho}+\Delta_{\sigma}U,$$ and, likewise:
$$tr_{g_0}D^2\varphi(x_0)=\beta^2\psi_{\rho
\rho}+\Delta_{\sigma}\psi.$$ For the linearized orthogonality
operator, note that, at $x_0$:
$${\cal L}_0{\cal O}[D\psi]=(\langle
\psi_{\rho},\tau_a\rangle+\psi_a,n\rangle )\tau_a.$$ Putting
everything together, the linear system to consider at $x_0$ is:
$$\left \{ \begin{array}{l}U_t-\beta^2U_{\rho
\rho}-\Delta_{\sigma}U=0,\\
\psi_t-\beta^2\psi_{\rho
\rho}-\Delta_{\sigma}\psi=0,\end{array}\right .$$ with boundary
conditions: $$\left \{ \begin{array}{l}U|_{\rho=0}=0,\\
\beta_0\langle \psi_{\rho},n\rangle+\beta
U_{\rho}|_{\rho=0}=b(\sigma,t),\\
\langle \psi_{\rho},\tau_a\rangle+\langle
\psi_a,n\rangle|_{\rho=0}=\omega_a(\sigma,t),\quad a=1,\ldots
n-1.\end{array} \right .$$ \vspace{.2cm}

Now take Fourier transform in $\sigma\in {\mathbb R}^{n-1}$
(corresponding to $\xi\in \mathbb{R}^{n-1}$), Laplace transform in
$t$ (corresponding to $p\in \mathbb{C}$) to obtain:
$$\hat{U}(\rho,\xi, p)\in \mathbb{C}, \hat{\psi}(\rho,\xi,p)\in
\mathbb{C}^n; \quad \xi\in{\mathbb R}^{n-1},p\in
\mathbb{C},\rho<0.$$ In transformed variables, we obtain the system
of linear ODE (in $\rho<0$, for fixed $(\xi,p)$):
$$\left \{ \begin{array}{l}\beta^2\hat{U}_{\rho
\rho}-(p+|\xi|^2)\hat{U}=0,\\
\beta^2\hat{\psi}_{\rho \rho}-(p+|\xi|^2)\hat{\psi}=0\end{array}
\right.$$ Writing the solution in the form:
 $$\left[\begin{array}{c}\hat{U}(\rho)\\\hat{\psi}(\rho)\end{array}\right]
 =e^{i\rho\gamma}\left[\begin{array}{c}\hat{U}(0)\\\hat{\psi}(0)\end{array}\right],$$
 we find the characteristic equation $\beta^2\gamma^2+p+|\xi|^2=0$,
 and choose the root $\gamma$ so that
 $i\gamma=(1/\beta)\sqrt{\Delta}$ (where $\Delta=p+|\xi|^2$ and we
 take the branch of $\sqrt{}$: $Re(\sqrt{\Delta})>0$). Here $(p,\xi)\in {\cal A}$, where:
 $${\cal A}=\{(p,\xi)\in {\mathbb C}\times \mathbb {R}^{n-1};
 |p|+|\xi|>0, Re(p)>-|\xi|^2\}.$$ Thus the solutions decay as
 $\rho\rightarrow -\infty$. Let ${\cal W}^+$ be the space of such
 decaying solutions, dim$_{\mathbb C}{\cal W}^+=n-1$. The relevant
 boundary operator on ${\cal W}^+$ is:
 $${\mathbb B}\left[\begin{array}{c}\hat{U}\\\hat{\psi}\end{array}\right]
 =\left[\begin{array}{c}\hat{U}\\ \beta_0\langle
 \hat{\psi}_{\rho},n\rangle+\beta \hat{U}_{\rho}\\
 \langle \hat{\psi}_{\rho},\tau_a\rangle+i\xi_a\langle
 \hat{\Psi},n\rangle\end{array}\right]_{|\rho=0}=
 \left[\begin{array}{c}\hat{U}(0)\\\beta_0(i\gamma)\langle
 \hat{\psi}(0),n\rangle +i\beta\gamma\hat{U}(0)\\
 (i\gamma)\langle \hat{\psi}(0),\tau_a\rangle+i\xi_a\langle
 \hat{\psi}(0),n\rangle \end{array}\right]$$(a vector in ${\mathbb
 C}\times{\mathbb C}\times{\mathbb C}^{n-1}).$

 The `complementarity condition' (see e.g. \cite{EidelmanZhitarasu}) is the statement that $\mathbb B$
 is a linear isomorphism from ${\cal W}^+$ to $\mathbb{C}^{n+1}$.
 With respect to the basis $\{\hat{U}(0),\langle
 \hat{\psi}(0),n\rangle,\langle \hat{\psi}(0),\tau_a\rangle\}$ of ${\cal
 W}^+$, the matrix of $\mathbb B$ is (in `block form'):
 $$[\mathbb{B}]=\left [\begin{array}{ccc}1&0&[0]_{1\times
 (n-1)}\\-\sqrt{\Delta}&-\frac{\beta_0}{\beta}\sqrt{\Delta}&[0]_{1\times
 (n-1)}\\ \left [0\right ]_{(n-1)\times 1}&[i\xi_a]_{(n-1)\times
 1}&-\frac{\sqrt{\Delta}}{\beta}\mathbb{I}_{n-1}\end{array}\right].$$
 This is triangular with non-zero diagonal entries for every
 $(p,\xi)\in {\cal A}$. Hence $\mathbb{B}$ is an isomorphism.

 \vspace{.4cm}
 \textbf{7. Estimates in H\"{o}lder spaces.}

 For the fixed-point argument based on the linear system, we need
 estimates for $||\cal F||_{\alpha}$, $||{\cal B}||_{1+\alpha}$,
 $||\Omega||_{1+\alpha}$, of two types: `mapping' and `contraction'
 estimates.\vspace{.2cm}

 A bit more precisely, for $T>0$, $R>0$ and $Q^T=D_0\times [0,T]$ consider the open ball:
 $$B_R^T=\{F\in C^{2+\alpha,1+\alpha/2}(Q^T,{\mathbb
R}^{n+1});||F-F_0||_{2+\alpha}<R, F|_{t=0}=F_0\}.$$
($F_0=[\varphi_0,u_0]$ is defined from the initial surface
$\Sigma_0$, via Lemma 4.1.) Solving the linear system with
`right-hand side' defined by $\bar{F}\in B_R^T$ defines a map
${\mathbb F}: \bar{F}\mapsto F$, and we need to verify that, for
suitable choices of $T$ and $R$, $\mathbb{F}$ maps into $B_R^T$ and
is a contraction.\vspace{.2cm}

\emph{Remark:} The argument that follows is standard, and the
experienced reader may want to skip to the statement of local
existence at the end of the next section. On the other hand, the
result is not covered by any general theorem proved in detail in a
reference known to the author, and some readers may find it useful
to have all the details included. Another reason is that, although
the `right hand sides' are clearly quadratic, without explicit
expressions one might run into trouble with compositions (which
behave poorly in H\"{o}lder spaces), or when appealing to `Taylor
remainder arguments' if the domain is not convex.\vspace{.2cm}

 For `mapping', we need estimates of the form:
 $$||{\cal F}(\bar{F},F_0)||_{\alpha}+||{\cal
 B}(D\bar{F},DF_0)||_{1+\alpha}+||\Omega(D\bar{\varphi},D\varphi_0)||_{1+\alpha}\mbox{
 decays as }T\rightarrow 0_+,$$
 and for `contraction':
 $$||{\cal F}(F^1,F^0)-{\cal F}(F^2,F^0)||_{\alpha}+||{\cal
 B}(DF^1,DF^2)||_{1+\alpha}+||\Omega(D\varphi^1,D\varphi^2)||_{1+\alpha}\leq
 \mu(T)||F^1-F^2||_{2+\alpha},$$
 where $\mu(T)\rightarrow 0$ as $T\rightarrow 0_+$.\vspace{.2cm}

 \emph{Notation:} The
 $(\alpha,\alpha/2)$ norms are taken on $Q^T$, the
 $(1+\alpha,(1+\alpha)/2)$ norms on $\partial D_0\times [0,T]$).
 Double bars without an index refer to the $(2+\alpha,1+\alpha/2)$
 norm, single bars to supremum norms over $Q^T$, and parabolic norms are
 indexed by their spatial regularity ($\alpha$ for
 ($\alpha,\alpha/2$), etc.) In general, we use brackets for
 H\"{o}lder-type difference quotients.
  \vspace{.2cm}

 We deal with the estimates for the `forcing term' $\cal F$ first. Consider the map $${\cal
 G}:Imm(\mathbb{R}^n,\mathbb{R}^{n+1})\rightarrow GL_n$$ which
 associates to the linear immersion $A$ the inverse matrix of
 $(\langle A_i,A_j\rangle)_{i=1}^n$, inner products of the rows of
 $A$. ${\cal G}$ is smooth, in particular locally Lipschitz in the
 space ${\cal W}$ of linear immersions. Hence if $F^1,F^2$ are maps
 $Q^T\rightarrow \mathbb{R}^{n+1}$, such that
 $DF^i\in C^{\alpha, \alpha/2}(Q^T)$ and $DF^i(z)\in K$ for all
 $z\in Q^T$, where $K\subset {\cal W}$ is a fixed compact set, we
 have the bound:
 $$||{\cal G}(DF^1)-{\cal G}(DF^2)||_{\alpha}\leq
 c_K||D(F^1-F^2)||_{\alpha}.$$
 In fact our maps $F^i$ are in $C^{2+\alpha,1+\alpha/2}$, so
 $DF^i\in C^{1+\alpha,\frac{1+\alpha}2}$. From this higher
 regularity we obtain the decay as $T\rightarrow 0_+$.
 Assuming $F^1|_{t=0}=F^2|_{t=0}$, we have:
 $$|D(F^1-F^2)|\leq
 [D(F^1-F^2)]_t^{(\frac{1+\alpha}2)}T^{\frac{1+\alpha}2}.$$

 Now recall the elementary fact: if $D\subset {\mathbb R}^n$ is a uniformly
 $C^1$ domain (not necessarily convex or bounded), and $f\in
 C^1(D)$ with $\alpha\in (0,1)$, we have: $[f]^{(\alpha)}\leq
 C_D||f||_{C^1}$. (Here `uniformly $C^1$' means $D$ can be covered by countably many
 balls of a fixed radius, which are domains of $C^1$
 manifold-with-boundary local charts for $D$, with uniform $C^1$
 bounds for the charts and their inverses. The constant $C_D$
 depends on those bounds.) Applying this to $DF$, where
 $F=F^1-F^2$ vanishes identically at $t=0$, and assuming $T<1$:

 $$[DF]_x^{(a)}\leq c(|DF|+|D^2F|)\leq
 c([DF]_t^{(\frac{1+\alpha}2)}T^{\frac{1+\alpha}2}+[D^2F]_t^{(\frac{\alpha}2)}T^{\alpha/2})\leq
 c||F||T^{\alpha/2},$$
 (where $c$ depends on $D_0$) and similarly for the oscillation in $t$:
 $$[DF]_t^{(\frac{\alpha}2)}\leq
 [DF]_t^{(\frac{1+\alpha}2)}T^{1/2}\leq ||F||T^{1/2},$$
 so we have:
 $$||D(F^1-F^2)||_{\alpha}\leq c ||F^1-F^2||T^{\alpha/2}.$$
 We conclude, under the
 assumption $F^1=F^2$ at $t=0$:
 $$||{\cal G}(DF^1)-{\cal G}(DF^2)||_{\alpha}\leq
 c_K||F^1-F^2||T^{\alpha/2}.$$
 In particular, applying this to $\bar{F}$ and $F_0$, we find:
 $$||({\cal G}(D\bar{F})-{\cal G}(DF_0))D^2\bar{F}||_{\alpha}\leq
 c_K||\bar{F}-F_0||T^{\alpha/2}||\bar{F}||,$$
 and for $F^1$ and $F^2$ coinciding at $t=0$:
 $$||({\cal G}(DF^1)-{\cal G}(DF^2))D^2F^1||_{\alpha}\leq
 c_K||F^1-F^2||T^{\alpha/2}||F^1||,$$
 as well as:
 $$||({\cal G}(DF^2)-{\cal G}(DF_0))(D^2F^1-D^2F^2)||_{\alpha}\leq
 c_K||F^2-F_0||T^{\alpha/2}||F^1-F^2||,$$
 so we have the mapping and contraction estimates for ${\cal
 F}(\bar{F},F_0)$ and ${\cal F}(F^1,F_0)-{\cal
 F}(F^2,F_0)$.

 \textbf{Lemma 7.1.} Assume $\bar{F},F_0,F^1,F^2$ are in
 $C^{2+\alpha,1+\alpha/2}(Q^T;{\mathbb R}^{n+1})$ and have the same
 initial values, and that
 $D\bar{F},DF_0,DF^1,DF^2$ all take
 values in the compact subset $K$ of $Imm(\mathbb{R}^n,\mathbb
 {R}^{n+1})$. Then:
 $$||{\cal F}(\bar{F},F_0)||_{\alpha}\leq
 c_K||\bar{F}-F_0||||\bar{F}||T^{\alpha/2},$$
 $$||{\cal F}(F^1,F_0)-{\cal F}(F^2,F_0)||_{\alpha}\leq
 c_K(||F^1||+||F^2-F_0||)T^{\alpha/2}||F^1-F^2||.$$ In particular,
 if $\bar{F}\in B_R^T$:
 $$||{\cal F}(\bar{F},F_0)||_{\alpha}\leq c_0RT^{\alpha/2}.$$
 If $\bar{F}^1,\bar{F}^2\in B_R^T$, we have:
$$||{\cal F}(\bar{F}^1,F_0)-{\cal F}(\bar{F}^2,F_0)||_{\alpha}\leq
c_0T^{\alpha/2}||\bar{F}^1-\bar{F}^2||.$$
 (The constant $c_0$ depends only on the data at $t=0$, and we assume
 $T<1$, $R<1$).\vspace{.2cm}

 Turning to the orthogonality boundary condition, first observe
 that:
 $$\Omega(D\varphi^1,D\varphi^2)=\langle
 D^T\varphi^1,d_n\varphi^1\rangle-\langle
 D^T\varphi^2,d_n\varphi^2\rangle-{\cal L}_0{\cal
 O}[D\varphi^1-D\varphi^2]$$
 $$=\langle
 D^T(\varphi^1-\varphi^2),d_n\varphi^1\rangle+\langle D^T\varphi^2,d_n(\varphi^1-\varphi^2)\rangle-\langle
 d_n(\varphi^1-\varphi^2),D^T\varphi_0\rangle-\langle
 D^T(\varphi^1-\varphi^2),d_n\varphi_0\rangle$$
 $$=\langle
 D^T\varphi^1-D^T\varphi^2,d_n\varphi^1-d_n\varphi_0\rangle
 +\langle
 d_n\varphi^1-d_n\varphi^2,D^T\varphi^2-D^T\varphi_0\rangle,$$
 which has quadratic structure. Using a local frame
 $(\tau_a)_{a=1}^{n-1}$ for $T\partial D_0$, we find the components
 $\Omega_a$:
 $$\Omega_a(D\varphi^1,D\varphi^2)=[\partial_i(\varphi^1-\varphi^2)\partial_j(\varphi^1-\varphi_0)
 +\partial_j(\varphi^1-\varphi^2)\partial_i(\varphi^2-\varphi_0)]n^j\tau_a^i$$
 (summation convention, $i,j=1,\ldots n$), so $\Omega_a$ is a sum of
 terms of the form:
 $b(x)D(\varphi^1-\varphi^2)D(\varphi^3-\varphi^4)$, where
 $b(x)=n^j\tau_a^i$ and the $\varphi^I$ coincide at $t=0$.
 It is then not hard to show that:
 $$||b(x)D(\varphi^1-\varphi^2)D(\varphi^3-\varphi^4)||_{1+\alpha}\leq
 c||b||_{1+\alpha}||\varphi^1-\varphi^2||||\varphi^3-\varphi^4||T^{\alpha},$$
 with $c$ depending on the $C^1$ norms of local charts for $D_0$.
 To bound the norm $||n\otimes\tau_a||_{1+\alpha}$, note $|n||\tau_a|\leq 1$,
 $|D(n\otimes \tau_a)|\leq |Dn|+|D\tau_a|$ and $[D(n\otimes
 \tau_a)]_{x}^{(\alpha)}\leq [Dn]_x^{(\alpha)}+[D\tau_a]_x^{\alpha}$.
 Since $n=-(\beta/\beta_0)Du_0$ on $\partial D_0$
 (and $\partial D_0$ is a level set of $u_0$), we clearly have:
 $$||Dn||_{\alpha}+||D\tau_a||_{\alpha}\leq c||D^2u_0||_{\alpha}\leq
 c||u_0||.$$ We summarize the conclusion in the following lemma.

\textbf{Lemma 7.2.} Assume
$\bar{\varphi},\varphi_0,\varphi^1,\varphi^2\in
C^{2+\alpha,1+\alpha/2}(Q^T;{\mathbb R}^n)$ have the same initial
values. Then:
 $$||\Omega(D\bar{\varphi},D\varphi_0)||_{1+\alpha}\leq
 c_0||u_0||||\bar{\varphi}-\varphi_0||^2T^{\alpha}$$
 and
 $$||\Omega(D\varphi^1,D\varphi^2)||_{1+\alpha}\leq
 c_0||u_0||(||\varphi^1-\varphi_0||+||\varphi^2-\varphi_0||)T^{\alpha}||\varphi^1-\varphi^2||,$$
 with $c_0$ depending only on the data at $t=0$. In particular, if
 $\bar{F}=[\bar{\varphi},\bar{u}]\in B_R^T$, we have:
$$||\Omega(D\bar{\varphi},D\varphi_0)||_{1+\alpha}\leq
c_0R^2T^{\alpha},$$ and for
$\bar{F}^I=[\bar{\varphi}^I,\bar{u}^I]\in B_R^T$, $I=1,2$:
$$||\Omega(D\bar{\varphi}^1,D\bar{\varphi}^2)||_{1+\alpha}\leq c_0
RT^{\alpha}||\bar{\varphi}^1-\bar{\varphi}^2||.$$
 \vspace{.3cm}

To explain the estimates for the angle condition, we write the
normal vector as a multilinear form on $DF^i$:
$$\tilde{N}(DF)=J_n(DF):=(-1)^n\sum_{i=1}^{n+1}(-1)^{i-1}(DF^1\wedge\ldots\hat{DF^i}\wedge\ldots
DF^{n+1})e_i\in{\mathbb R}^{n+1}$$ ($DF^i$ omitted in the $i^{th.}$
term of the sum), where $DF^i\in {\mathbb R}^n$ for $i=1,\ldots,
n+1$ and we identify the $n$-multivector in ${\mathbb R}^n$ with a
scalar, using the standard volume form. The angle condition has the
form:
$$\beta^2|\tilde N|^2-\langle \tilde{N},e_{n+1}\rangle^2=0\mbox{ on
}\partial D_0,$$ and we set:
$$B(DF):=\beta^2|J_n(DF)|^2-\langle J_n(DF),e_{n+1}\rangle^2,$$
with linearization at $DF_0=[\mathbb{I}_n|Du_0]$:
$${\cal L}_0B[DF]=2\beta^2\langle
J_n(DF_0),DJ_n(DF_0)[DF]\rangle-2\langle
J_n(DF_0),e_{n+1}\rangle\langle DJ_n(DF_0)[DF],e_{n+1}\rangle.$$

Under the assumption $F^1=F^2$ at $t=0$, we need an estimate in
$C^{1+\alpha,\frac{1+\alpha}2}$ for:
$${\cal B}(DF^1,DF^2):=B(DF^1)-B(DF^2)-{\cal L}_0B[DF^1-DF^2]$$
$$=\beta^2(|J_n(DF^1)|^2-|J_n(DF^2)|^2-2\langle
J_n(DF_0),DJ_n(DF_0)[DF^1-DF^2]\rangle)$$
$$-(\langle J_n(DF^1),e_{n+1}\rangle^2-\langle
J_n(DF^2),e_{n+1}\rangle^2-2\langle J_n(DF_0),e_{n+1}\rangle \langle
DJ_n(DF_0)[DF^1-DF^2],e_{n+1}\rangle ).$$ It will suffice to
estimate the expression in the first parenthesis; the second is
analogous.\vspace{.2cm}

We need the following algebraic observation: if
$T_0=[\mathbb{I}_n|Du_0]$ and $T$ are $n\times (n+1)$ matrices, the
expression:
$$|J_n(T_0+T)|^2-|J_n(T_0)|^2-2\langle
J_n(T_0),DJ_n(T_0)[T]\rangle$$ is a linear combination (with
constant coefficients) of terms of the form:
$$u_{0i}p_{(2)}(T),\quad u_{0i}u_{0j}p_{(2)}(T),\quad p_{(2)}(T),$$
where the $p_{(2)}(T)$ are polynomials in the entries of $T$ (with
constant coefficients), with terms of degree: $2\leq deg\leq 2n$.

Thus ${\cal B}(DF^1,DF^2)$ is a linear combination (with constant
coefficients) of terms:
$$u_{0i}p_{(2)}(DF^1-DF^2),\quad u_{0i}u_{0j}p_{(2)}(DF^1-DF^2),\quad p_{(2)}(DF^1-DF^2),$$
with the $p_{(2)}$ as described; and hence is a linear combination
of terms of the form:
$$u_{0i}(F^{1j}_k-F^{2j}_k)^d,\quad
u_{0i}u_{0l}(F^{1j}_k-F^{2j}_k)^d,\quad (F^{1j}_k-F^{2j}_k)^d$$
(where $2\leq d\leq 2n,1\leq j\leq n+1,1\leq i,l,k\leq n$), which we
write symbolically as:
$${\cal B}(DF^1,DF^2)\sim \sum_{2\leq d\leq2n}b(x)(DF^1-DF^2)^d,$$
where $b(x)$ is constant or $u_{0i}(x)$ or $u_{0i}(x)u_{0j}(x).$ For
the degree $d$ terms $G^{(d)}\sim b(x)(DF^1-DF^2)^d$, it is not hard
to show the bound:
$$||G^{(d)}||_{1+\alpha}\leq
c||b||_{1+\alpha}||F^1-F^2||^dT^{\alpha},\quad 2\leq d\leq 2n.$$ We
conclude:

\textbf{Lemma 7.3.}Assume $\bar{F},F_0,F^1,F^2$ are in
 $C^{2+\alpha,1+\alpha/2}(Q^T;{\mathbb R}^{n+1})$ and have the same
 initial values. Then:
$$||{\cal B}(D\bar{F},DF_0)||_{1+\alpha}\leq
c(1+||u_0||^2)(1+||\bar{F}-F_0||^{2n-2})T^{\alpha}||\bar{F}-F_0||^2.$$
$$||{\cal B}(DF^1,DF^2)||_{1+\alpha}\leq
c(1+||u_0||^2)(1+||F^1-F^2||^{2n-2})T^{\alpha}||F^1-F^2||^2,$$ with
$c$ depending only on $F_0$. In particular, if $\bar{F}\in B_R^T$:
$$||{\cal B}(D\bar{F},DF_0)||_{1+\alpha}\leq c_0 R^2T^{\alpha},$$
and if $\bar{F}^1,\bar{F}^2\in B_R^T$:
$$||{\cal B}(D\bar{F}^1,D\bar{F}^2)||_{1+\alpha}\leq
c_0T^{\alpha}||\bar{F}^1-\bar{F}^2||,$$ with $c_0$ depending only on
$F_0$.
 \vspace{.3cm}

\textbf{8.} \textbf{Local existence.}

Given a $C^{3+\bar{\alpha}}$ graph $\Sigma_0$ over a uniformly
$C^{3+\bar{\alpha}}$ domain $D_0\subset \mathbb{R}^n$ (for arbitrary
$\bar{\alpha}\in (0,1)$) satisfying the contact and angle
conditions, let $\varphi_0\in Diff^{2+\bar{\alpha}}$ be the
diffeomorphism given by lemma 4.1 (with the 1-jet of the identity at
$\partial D_0$ and 2-jet determined by the mean curvature of
$\Sigma_0$ at $\partial D_0$). Then find $u_0\in C^{2+\alpha}(D_0)$
so that $F_0=[\varphi_0,u_0]\in C^{2+\alpha}(D_0;\mathbb{R}^{n+1})$
parametrizes $\Sigma_0$ over $D_0$
($\alpha=\bar{\alpha}^2<\bar{\alpha}$).

(Precisely, if $[z,\tilde{u_0}(z)]$ parametrizes $\Sigma_0$ as a
graph, and $\varphi_0$ is given by lemma 4.1, let
$u_0=\tilde{u_0}\circ \varphi_0$; so $u_0\in C^{2+\alpha}$.)
 \vspace{.2cm}

We obtained in section 7 all the estimates needed for a fixed-point
argument in the set:
$$B_R^T=\{F\in C^{2+\alpha,1+\alpha/2}(Q^T,{\mathbb
R}^{n+1});||F-F_0||<R, F|_{t=0}=F_0\}.$$

Choose $R<1$ and $T_0<1$ small enough (depending only on $F_0$) so
that, for $F\in B_R^{T_0}$, $F(t)=[\varphi(t),u(t)]$ defines an
embedding of $D_0$, with $\varphi(t)$ a diffeomorphism onto its
image $D(t)$. Let $K\subset Imm(\mathbb{R}^n,\mathbb
 {R}^{n+1})$ be a compact set containing $DF(z)$ for all $F\in B_R,
 z\in Q^{T_0}$. Now consider $T<T_0$.
\vspace{.2cm}

Given $\bar{F}\in B_R^T$, solve the linear system (LPDE/LBC) (with
initial data $F_0$) to obtain $F\in C^{2+\alpha,1+\alpha/2}(Q^T)$.
(This is possible since the complementarity and compatibility
conditions hold for the linear system.) This defines a map ${\mathbb
F}:{\bar F}\mapsto F$.

From linear parabolic theory (e.g. \cite{EidelmanZhitarasu}, thm
VI.21]):
$$||F-F_0||\leq M(||{\cal F}(\bar{F},F_0)||_{\alpha}+||{\cal
B}(D\bar{F},DF_0)||_{1+\alpha}+||\Omega(D\bar
{\varphi},D\varphi_0)||_{1+\alpha}),$$ where $M>0$ depends on the
$C^{\alpha,\alpha/2}$ norm of the coefficients of the linear system,
that is, ultimately on $||F_0||$. \vspace{.2cm}

From lemmas 7.2-7.4 in section 7, it follows that:
$$||F-F_0||\leq Mc_0(RT^{\alpha/2}+R^2T^{\alpha})<R$$
provided $T$ is chosen small enough (depending only on $F_0$.) Thus
$\mathbb{F}$ maps $B_R^T$ to itself.\vspace{.2cm}

Similarly, if ${\mathbb F}(\bar{F}^i)=F^i$ for $i=1,2$, standard
estimates for the linear system solved by $F^1-F^2$ give:
$$||F^1-F^2||\leq M(||{\cal F}(\bar{F}^1,\bar{F}^2)||_{\alpha}+||{\cal
B}(D\bar{F}^1,D\bar{F}^2)||_{1+\alpha}+||\Omega(D
{\bar{\varphi}}^1,D\bar{\varphi}^2)||_{1+\alpha})$$

Again the estimates in lemmas 7.2-7.4 imply:
$$||F^1-F^2||\leq Mc_0(T^{\alpha/2}+T^{\alpha})||\bar{F}^1-\bar{F}^2||<\frac 12||\bar{F}^1-\bar{F}^2||,$$
assuming $T$ is small enough (depending only on $F_0$). This
concludes the argument for local existence. \vspace{.2cm}

\textbf{Theorem 8.1.} Let $\Sigma_0\subset {\mathbb R}^{n+1}$ be a
$C^{3+\bar{\alpha}}$ graph over $D_0\subset {\mathbb R}^n$
satisfying the contact and angle conditions at $\partial D_0$. With
$\alpha=\bar{\alpha}^2$, there exists a parametrization
$F_0=[\varphi_0,u_0]\in C^{2+\alpha}(D_0)$ of $\Sigma_0$, a number
$T>0$ depending only on $F_0$ and a unique solution $F\in
C^{2+\alpha,1+\alpha/2}(Q^T;{\mathbb R}^{n+1})$ of the system:
$$\left \{ \begin{array}{l}
\partial_t F-g^{ij}(DF)\partial_i\partial_jF=0,\quad F=[\varphi,u]\\
u_{|\partial D_0}=0, \quad N^{n+1}(D\varphi,Du)_{|\partial
D_0}=\beta,\quad \langle D^T\varphi,d_n\varphi\rangle_{|\partial
D_0}=0\end{array}\right .$$ with initial data $F_0$. For each $t\in
[0,T]$, $F(t)$ is a $C^{2+\alpha}$ embedding parametrizing a surface
$\Sigma_t$ which satisfies the contact and angle conditions and
moves by mean curvature. In addition, $F(t)$ satisfies the
orthogonality condition at $\partial D_0$.

The hypersurfaces $\Sigma_t$ are graphs. For each $t\in [0,T]$,
$\varphi(t):D_0\rightarrow D(t)$ is a diffeomorphism and
$\Sigma_t=graph(w(t))$, for $w(t):D(t)\rightarrow {\mathbb R}$ given
by $w(t)=u(t)\circ \varphi^{-1}(t)$. (We have $w(t)\in
C^{2+\alpha^2}(D(t))$, `less regular' than $u(t)$ or $\varphi(t)$.)
$D(t)$ is a uniformly $C^{2+\alpha}$ domain.\vspace{.3cm}

\emph{Remark 8.1.} This theorem does not address geometric
uniqueness of the motion, given $\Sigma_0$. It only asserts
uniqueness for solutions of the parametrized flow (including the
orthogonality boundary condition) in the given regularity
class.\vspace{.3cm}

\textbf{9. Rotational symmetry.} In this section we record the
equations for two rotationally symmetric instances of the problem:
(i) $D_0$ and $D(t)$ are disks, and $u>0$ (`lens' case); (ii) $D_0$
and $D(t)$ are complements of disks in ${\mathbb R}^n$ (`exterior'
case). For simplicity we restrict to $n=2$.

Let $F(r)=[\varphi(r),u(r)]$ parametrize a hypersurface $\Sigma$,
where $\varphi(r)=\phi(r)e_r$ is a diffeomorphism onto its image.
Here $e_r,e_{\theta}$ are orthonormal vectors, outward normal (resp.
counterclockwise tangent) to the circles $r$=const. The unit upward
normal vector and mean curvature are:
$$N=\frac{[-u_re_r,\phi_r]}{\sqrt{u_r^2+\phi_r^2}},$$
$$H=\frac 1{(\phi_r^2+u_r^2)^{3/2}}(\phi_r{\cal
M}(\phi_r,u_r)[D^2u]-\langle u_re_r,\vec{\cal
M}(\phi_r,u_r[D^2\varphi]\rangle),$$ where:
$${\cal
M}(\phi_r,u_r)[D^2u]=u_{rr}+(\phi_r^2+u_r^2)\frac{u_r\phi_r}{\phi^2},$$
$$\vec{\cal
M}(\phi_r,u_r)[D^2\varphi]=[\phi_{rr}+(\phi_r^2+u_r^2)(\frac{r\phi_r}{\phi^2}-\frac
1{\phi})]e_r.$$ Simplifying:
$$H=\frac
1{(\phi_r^2+u_r^2)^{3/2}}[\phi_ru_{rr}-u_r\phi_{rr}+(\phi_r^2+u_r^2)\frac{u_r}{\phi}].$$

Now consider the time-dependent case $F(r,t)=[\phi(r,t)e_r,u(r,t)]$.
From the above expressions, one finds easily that the equation
$\langle\partial_tF,N\rangle=H$ takes the form:
$$\phi_r(u_t-\frac 1{\phi_r^2+u_r^2}{\cal
M}(\phi_r,u_r)[D^2u])=u_r\langle e_r,\varphi_t-\frac
1{\phi_r^2+u_r^2}\vec{\cal M}(\phi_r,u_r)[D^2\varphi]\rangle.$$

In `split gauge', we consider the system:
$$\left \{ \begin{array}{l}u_t-\frac 1{\phi_r^2+u_r^2}{\cal
M}(\phi_r,u_r)[D^2u]=0,\\
\varphi_t-\frac 1{\phi_r^2+u_r^2}\vec{\cal
M}(\phi_r,u_r)[D^2\varphi]=0.\end{array}\right .$$

Note that $\phi(r,t)=r$ solves the $\phi$ equation, and that in this
case the $u$ equation becomes:
$$w_t-\frac{w_{rr}}{1+w_r^2}-\frac{w_r}r=0.$$ This can be compared
with the equation for curve networks:
$$w_t-\frac{w_{xx}}{1+w_x^2}=0.$$\vspace{.2cm}

 The boundary conditions are easily stated (we assume $D_0$ is the unit disk or its
 complement).

 The `contact condition' at $r=1$
is $u=0$. For the `angle condition' at $r=1$, we find:
$$u_r^2=\frac{\beta_0^2}{\beta^2}\phi_r^2,\quad
\beta_0:=\sqrt{1-\beta^2}.$$ Assuming $\phi_r>0$ at $r=1$, this
resolves as:
$$\beta u_r+\beta_0\phi_r=0\mbox{ at }r=1\mbox{ (lens case)};$$
$$\beta u_r-\beta_0\phi_r=0\mbox{ at }r=1\mbox{ (exterior case)}.$$
(For lenses, one also has at $r=0$: $u_r=0$ and $\phi_r=1$.) Thus in
both cases one can work with \emph{linear} Dirichlet/Neumann-type
boundary conditions.\vspace{.2cm}

One reason to consider the exterior case is that (unlike the lens
case) it admits stationary solutions. Geometrically, one just has to
consider one-half of a catenoid, truncated at an appropriate height.
For example, for 120 degree junctions the equation for stationary
solutions:
$$\left \{ \begin{array}{l}\frac{u_{rr}}{1+u_r^2}+\frac{u_r}r=0\mbox{ in
}\{r>1\},\\
{u_r}_{|r=1}=\sqrt{3},\quad u_{|r=1}=0\end{array}\right .$$ admits
the explicit solution:
$$u(r)=\frac{\sqrt{3}}2(\ln(2r+\sqrt{4r^2-3})-\ln 3),\quad
r>\sqrt{3}/2.$$\vspace{.3cm}

\emph{Problem.} It would be interesting to consider the nonlinear
dynamical stability of this solution (even linear stability is yet
to be considered.) One may even work with bounded domains, by
introducing a fixed boundary at some $R>1$, intersecting the surface
orthogonally (see Section 10.)\vspace{.3cm}

\textbf{10. Fixed supporting hypersurfaces.} Extending the local
existence theorem to the case of hypersurfaces intersecting a fixed
hypersurface $\cal S$ orthogonally presents no essential difficulty.
The case of vertical support surface leads directly to graph
evolution with a standard Neumann condition on a fixed boundary; we
consider the complementary case where $\cal S$ is a graph. Let
${\cal S}\subset {\mathbb R}^{n+1}$ be a $C^{4}$ embedded
hypersurface (not necessarily connected), the graph over ${\cal
D}\subset {\mathbb R}^n$ of $B\in C^4({\cal D})$, oriented by the
upward unit normal:
$$\nu(y):=\frac 1{v_B}\tilde{\nu}(y),\quad
\tilde{\nu}(y):=[-DB(y),1]\in {\mathbb R}^n\times {\mathbb R}, \quad
v_B:=\sqrt{1+|DB(y)|^2}.$$ $\nu$ is assumed to be nowhere vertical
in ${\cal D}$ ($DB\neq 0$).  To state the problem in the graph
parametrization, we consider a time-dependent domain $D(t)\subset
{\mathbb R}^n$ with boundary consisting of two components
$\partial_1D(t)$ and $\partial_2D(t)$, both moving. The hypersurface
$\Sigma_t$ is the graph of $w(\cdot, t)$ over $D(t)$, solving the
parabolic equation:
$$w_t-g^{ij}(Dw)w_{ij}=0\mbox{ in }E:=\bigcup_{t\in [0,T]}D(t)\times
\{t\}\in {\mathbb R}^{n+1}\times [0,T],$$ with boundary conditions:
$$w(\cdot,t)_{|\partial_1D(t)}=0,\quad
\sqrt{1+|Dw|^2}_{|\partial_1D(t)}=1/\beta$$ (as before), and on
$\partial_2D(t)$:
$$w=B, \quad \nabla w\cdot \nabla B=-1.$$
 (The first-order condition on $\partial_2D(t)$ is equivalent to
$\langle \nu,N\rangle=0$).\vspace{.2cm}

Differentiating in $t$ the boundary condition $w=B$ leads easily to
an equation for the normal velocity of the interface
$\Gamma(t)=\partial_2D(t)$:
$$\dot{\Gamma_n}=\frac{vH}{B_n-w_n}.$$\vspace{.2cm}
Note that $w_n$ at $\partial_2D(t)$ can be computed from $B_n$,
since:
$$-1=\nabla w\cdot \nabla B=w_nB_n+|\nabla^TB|^2;$$
in particular neither $B_n$ nor $w_n$ can vanish (so both have
constant sign on connected components of $\partial_2D$), and one
easily computes: $w_n-B_n=-v_B^2/B_n$.\vspace{.2cm}

Let $\Lambda=\Sigma\cap {\cal S}$ be the $(n-1)$-manifold of
intersection, the graph of $w$ (or $B$) over $\partial_2D$. Given
the graph parametrizations of $\Sigma$ and $\cal S$:
$$G(y)=[y,w(y)],\quad {\mathbb B}(y)=[y,B(y)],\quad y\in
\partial_2D,$$
and $\tau \in T\partial_2 D$, we have the tangent vectors:
$$G_n:=[n,w_n]\in T\Sigma,\quad G_B:=[\nabla B,-1]=-v_B\nu\in
T\Sigma,\quad G_{\tau}:=[\tau, \nabla w\cdot \tau]\in T\Lambda,$$
and the second fundamental forms of $\Sigma$ and ${\cal S}$(for
$e\in {\mathbb R}^n$ arbitrary):
$$A(dGe,dGe)=\frac 1{v}d^2w(e,e),\quad {\cal A}(d{\mathbb
B}e,d{\mathbb B}e)=\frac 1{v_B}d^2B(e,e).$$ From $\langle
\nu,N\rangle=0$ at $\partial_2D$, it follows easily that (cp.
\cite{Stahl}):
$$A(G_\tau,\nu)=-{\cal A}(G_\tau,N),
\quad \tau\in T\partial D.$$

For the remainder of this section, we concentrate on the boundary
conditions at $\partial_2D_0$. To establish short-time existence, we
consider as before the parametrized flow:
$$F_t-tr_gd^2F=0,\quad g=g(dF),\quad F=[\varphi,u].$$
The contact and angle boundary conditions are:
$$u_{|\partial_2 D_0}=B\circ \varphi_{|\partial_2 D_0},\quad \langle
N,\nu\circ \varphi\rangle_{|\partial_2 D_0}=0.$$ Again we have two
scalar boundary conditions for $n+1$ components. Here the solution
is easier than at the junction. With the notation
$F_n=dFn=[\varphi_n,u_n]$, we replace the angle condition by the
`vector Neumann condition':
$$F_n\perp T{\cal S},\mbox{ or }F_n=-\alpha v_B\nu\mbox{ on
}\partial_2 D_0,$$ where $\alpha:\partial_2 D_0\rightarrow {\mathbb
R}$, or equivalently (since this leads to $\alpha=-u_n$):
$$\varphi_n=-u_n(\nabla B\circ \varphi) \mbox{ on }\partial_2D_0.$$
Clearly the Neumann condition implies the angle condition $\langle
N,\nu\circ \varphi\rangle=0$, but not conversely. This linear
Neumann-type condition can easily be incorporated into the
fixed-point existence scheme described earlier. \vspace{.2cm}

There is one issue to consider: the 0 and 1st-order compatibility
conditions must hold at $\partial_2 D_0$, at $t=0$. The initial
hypersurface $\Sigma_0$ uniquely determines $w_0$ and $D_0\subset
\mathbb{R}^n$ (satisfying $w_0=B$ and $\nabla w_0\cdot \nabla B=-1$
on $\partial_2 D_0$), and then once $\varphi_0\in Diff(D_0)$ is
fixed, $u_0=w_0\circ \varphi_0$ is also determined. We may assume:
$${\varphi_0}=id, \quad \varphi_{0n}=\nabla B\mbox{ on }\partial_2
D_0,$$ so:$$u_{0n}=\nabla w_0\cdot \varphi_{0n}=\nabla w_0\cdot
\nabla B=-1\mbox{ on }\partial_2 D_0,$$ and then the Neumann
condition ${F_{0n}}_{|\partial_2D_0}=-v_B\nu$ holds at $t=0$, on
$\partial_2 D_0$.\vspace{.2cm}

The first-order compatibility condition is:
$$tr_gd^2u_0=u_t=\nabla B\cdot \varphi_t=\nabla B\cdot
tr_gd^2\varphi_0\mbox{ on }\partial D_0,$$ or equivalently:
$$tr_g\langle \nu,d^2F_0\rangle=0\mbox{ on }\partial D_0.$$
(This is not a mean curvature condition; the mean curvature of
$\Sigma_0$ is $H=tr_g\langle N,d^2F_0\rangle.$) \vspace{.2cm}

From now on we omit the subscript $0$, but continue to discuss
compatibility at $t=0$. First observe that the Neumann condition
leads to a splitting of the induced metric. Given $\tau \in
T\partial_2 D_0$, let $F_{\tau}=dF\tau\in T\Lambda$. Then (recalling
$u_n=-1$ on $\partial_2 D_0$):
$$\langle F_{\tau},F_n\rangle=\langle[\tau, dB
\tau],[\varphi_n,u_n]\rangle=\nabla B\cdot \tau-\nabla B\cdot
\tau=0.$$ Thus we have:
$$tr_g\langle \nu,d^2F\rangle=g^{ab}\langle
\nu,d^2F(\tau_a,\tau_b)\rangle +g^{nn}\langle
\nu,d^2F(F_n,F_n)\rangle,$$ for a local basis
$\{T_a=dF\tau_a\}_{a=1}^{n-1}$ of $T\Lambda$, with $g_{ab}=\langle
T_a,T_b\rangle$ and $g_{nn}=|F_n|^2=v_B^2$.\vspace{.2cm}

Differentiating in $n$ the condition $u_n=\nabla w\cdot \varphi_n$ (
assuming, as usual, $n$ extended to a tubular neighborhood $\cal N$
of ${\partial_2 D}_0$ as a self-parallel vector field), we find:
$$u_{nn}=d^2w(n,\nabla B)+ \nabla w\cdot d^2\varphi(n,n).$$
This is used to compute:
$$\langle \nu,d^2F(n,n)\rangle =\frac 1{v_B}[u_{nn}-\nabla B\cdot
d^2\varphi(n,n)]$$
$$=\frac 1{v_B}[d^2w(n,\nabla B)+(\nabla w-\nabla
B)\cdot d^2\varphi(n,n)]$$
$$=-vA(G_n,\nu)+\frac 1{v_B}(w_n-B_n)
n\cdot d^2\varphi(n,n).$$ Bearing in mind the expression for
$w_n-B_n$ found earlier, the compatibility condition may be stated
in the form:
$$\frac{v_B}{B_n}
n \cdot d^2\varphi(n,n)=-vA(G_n,\nu)+g^{ab}\langle
d^2F(\tau_a,\tau_b),\nu\rangle.$$\vspace{.2cm}

We are now in the same situation as in section 4: given the 1-jet of
$\varphi_0$ on $\partial_2 D_0$, we extend $\varphi_0$ to a tubular
neighborhood $\cal N$ of $\partial_2 D_0$ (and then to all of
$D_0$), so that $ n\cdot d^2\varphi(n,n)$ has on $\partial_2 D_0$
the value dictated by the compatibility condition (using Lemma
4.1(ii)). We just need to verify that the right-hand side of the
above expression depends only on $\Sigma_0$, $\cal S$, and the 1-jet
of $\varphi_0$ over $\partial_2 D_0$. Clearly only the term
$g^{ab}\langle \nu,d^2F(\tau_a,\tau_b)\rangle$ is potentially an
issue. \vspace{.2cm}

Fix $p\in \partial_2 D_0$, and let $\{\tau_a\}$ be an orthonormal
frame for $T\partial_2 D_0$ near $p$, parallel at $p$ for the
connection induced on $\partial_2 D_0$ from ${\mathbb R}^n$. If
$\cal K$ denotes the second fundamental form of $\partial_2 D_0$ in
${\mathbb R}^n$, we have:
$$\tau_a(\tau_b)={\cal K}(\tau_a,\tau_b)n\quad (\mbox{ at }p)$$
(on the left-hand-side, $\tau_b$ is regarded as a vector-valued
function in ${\mathbb R}^n$). Still computing at $p$, this implies:
$$d^2F(\tau_a,\tau_b)=\tau_a(dF\tau_b)-dF(\tau_a(\tau_b))$$
$$=\tau_a(d{\mathbb B}\tau_b)-{\cal K}(\tau_a,\tau_b)F_n$$
$$=d^2{\mathbb B}(\tau_a,\tau_b)+{\cal K}(\tau_a,\tau_b){\mathbb
B}_n-{\cal K}(\tau_a,\tau_n)F_n,$$ where $F_n=-v\nu$ and ${\mathbb
B}_n=d{\mathbb B}n\in T{\cal S}$. Hence:
$$\langle \nu,d^2F(\tau_a,\tau_b)\rangle=\langle \nu,d^2{\mathbb
B}(\tau_a,\tau_b)\rangle +v{\cal K}(\tau_a,\tau_b)={\cal
A}(T_a,T_b)+v{\cal K}(\tau_a,\tau_b).$$ This clearly depends only on
$\cal S$ and on $\Sigma_0$. \vspace{.2cm} We summarize the
discussion in a lemma.\vspace{.2cm}

\textbf{Lemma 10.1} Let $\Sigma_0=graph(w_0)$ be a $C^3$ graph over
$D_0\subset {\mathbb R}^n$ (a uniformly $C^3$ domain), intersecting
a fixed hypersurface ${\cal S}=graph(B)$ over $\partial D_0$.
Consider the parametrized mean curvature motion with Neumann
boundary condition:
$$F\in C^{2,1}(D_0\times [0,T])\rightarrow {\mathbb R}^{n+1},\quad
F=[\varphi,u]$$
$$F_t-tr_gd^2F=0,\quad g=g(dF),\quad u\circ \varphi=B\mbox{ and }F_n\perp T{\cal S}\mbox{
on }\partial D_0.$$ Then $\varphi_0\in Diff(D_0)$ can be chosen so
that (with $u_0=w_0\circ \varphi_0$) the initial data
$F_0=[\varphi_0,u_0]$ satisfies the order zero and the first-order
compatibility conditions at $t=0$ and $\partial D_0$:
$$\varphi_{0n}=-u_{0n}(\nabla B\circ {\varphi_0}),\quad \langle
\nu\circ \varphi_0,tr_{g_0}d^2F_0\rangle=0.$$

\emph{Remark 10.1.} Differentiating $dw\tau_a=dB\tau_a$ along
$\tau_b$, we find:
$$d^2w(\tau_a,\tau_b)-d^2B(\tau_a,\tau_b)=(w_n-B_n){\cal
K}(\tau_a,\tau_b)$$ (reminding us that, although $w\equiv B$ on
$\partial D_0$, the tangential components of their Hessians do not
coincide.) From this follows the expression for ${\cal K}$ in terms
of $A$ and $\cal A$:
$${\cal K}(\tau_a,\tau_b)=\frac
1{w_n-B_n}[vA(T_a,T_b)-v_B{\cal A}(T_a,T_b)].$$ It is also easy to
express the corresponding traces in terms of the mean curvatures
$H^{\Lambda}$ and ${\cal H}^{\Lambda}$ of $\Lambda$ in $\Sigma$ and
${\cal S}$:
$$H^{\Lambda}=\frac v{v_B}g^{ab}A(T_a,T_b),\quad {\cal
H}^{\Lambda}=\frac{v_B}vg^{ab}{\cal A}(T_a,T_b).$$\newpage

\textbf{11. Boundary conditions for the second fundamental form.}
\vspace{.2cm}

To understand the long-term behavior of a graph $(\Sigma_t)$ in
$\mathbb{R}^{n+1}$ moving by mean curvature and intersecting
$\mathbb{R}^n$ at a constant angle, we need to consider the
evolution of its second fundamental form. Working in the graph
parametrization the boundary conditions are easy to state (and
linear):
$$w_{|\partial D(t)}=0,\qquad d_nw_{|\partial
D(t)}=\frac{\beta_0}{\beta},$$ where $n=n_t$ is the inner unit
normal to $\partial D(t)$. It is possible to reparametrize the
$\Sigma_t$ over a different time-dependent domain ${\cal D}(t)$,
obtaining mean curvature flow:
$${\cal F}_t:{\cal D}(t)\rightarrow \mathbb{R}^{n+1},\qquad
\partial_t{\cal F}=HN,$$
with boundary conditions:
$${\cal F}^{n+1}_{|\partial {\cal D}(t)}=0,\qquad N^{n+1}_{|\partial {\cal
D}(t)}=\beta.$$ For this parametrization the evolution equation for
the second fundamental form (and its covariant derivatives of
arbitrary order) is well-understood \cite{Huisken}. The disadvantage
is that the unit normal $N_{|{\cal F}_t}$ depends non-linearly on
the components of $\cal F$, and as a result the boundary conditions
for the second fundamental form (which are needed for global
estimates over spacetime domains) do not admit simple expressions.
Therefore we choose to work with graph flow, at the cost of having
to derive and understand a new set of evolution equations. The
equations for $h$ and the mean curvature $H$ are derived in an
appendix. In this section we derive boundary conditions; the
development is similar to work of A. Stahl \cite{Stahl} for MCF of
hypersurfaces intersecting a fixed boundary
orthogonally.\vspace{.2cm}

It is easy to see that $h$ splits on $\partial D(t)$: if $\tau \in
T\partial D(t)$ is a tangential vector field, and $n=n_t$ is the
inner unit normal:
$$h(n,\tau)=\frac 1vd^2w(n,\tau)=\frac 1v(\tau(w_n)-Dw\cdot
\bar{\nabla}_\tau n)=0\mbox{ on }\partial D(t),$$ since $w_n\equiv
\beta_0/\beta$ on the boundary and $\bar{\nabla}_{\tau}n\in
T\partial D(t)$ ($\bar{\nabla}$ is the euclidean connection.) In
particular, it follows that $h(Dw,\tau)=0$ on $\partial D(t)$.
\vspace{.2cm}

\emph{Remark.} Already this simple fact cannot be shown for
$a(\nu,\tau)$, the second fundamental form in the MCF
parametrization, regarded as a quadratic form on ${\cal
D}(t)$.\vspace{.2cm}

\emph{Boundary condition for H.} We derived in section 2 the
equation for the normal velocity of the moving boundary
$\Gamma_t=\partial D(t)$:
$$\dot{\Gamma}_n=-\frac v{w_n}H=-\frac 1{\beta_0}H\quad \mbox{ at
}\partial D(t).$$ Since $\langle N,e_{n+1}\rangle
(\Gamma(t),t)\equiv \beta$ on $\partial D(t)$ we have: $$\langle
\partial_tN,e_{n+1}\rangle=-\langle \partial_kN,e_{n+1}\rangle\dot{\Gamma}^k,$$
where $\partial_kN=-g^{ij}h_{ik}G_j$, with $e_{n+1}$ component:
$$\langle \partial_kN,e_{n+1}\rangle=-g^{ij}w_jh_{ik}=-\frac
1{v^2}h(Dw,\partial_k)=-\frac 1{v^2}w_nh(n,\partial_k).$$ Hence we
find, on $\partial D(t)$:
$$\langle \partial_tN,e_{n+1}\rangle=\frac
{w_n}{v^2}h(n,\dot{\Gamma})=\frac
{w_n}{v^2}\dot{\Gamma}_nh(n,n)=-\beta Hh_{nn}.\qquad (11.1)$$ (We
set $h_{nn}:=h(n,n)$). Denote by $\nabla^{\Sigma}$ the gradient of
$\Sigma_t$, in the induced metric ($\nabla^{\Sigma}f=g^{ij}f_iG_j$).
Using $\partial_tN=-\nabla^{\Sigma}H-Hv^{-1}\nabla^{\Sigma }v$,
combined with the expressions (valid on $\partial D(t)$):
$$\langle \nabla^{\Sigma}H,e_{n+1}\rangle=g^{ij}{H_i}\langle
G_j,e_{n+1}\rangle=g^{ij}H_iw_j=\frac
1{v^2}w_iH_i=\frac{w_n}{v^2}H_n=\beta \beta_0H_n,$$
$$\langle
\nabla^{\Sigma}v,e_{n+1}\rangle=\frac{v_nw_n}{v^2}
=\frac{w_n^2}{v^2}h_{nn}=\beta_0^2h_{nn},$$ we find on $\partial
D(t)$:
$$\langle
\partial_tN,e_{n+1}\rangle=-\beta\beta_0(H_n+\beta_0Hh_{nn}).\qquad (11.2)$$
Comparing expressions for (11.1) and (11.2) for $\langle
\partial_tN,e_{n+1}\rangle$ yields a Neumann-type condition for $H$.
We state this as a lemma (including the evolution equation derived
in the appendix). Here $L=L_g$ denotes the operator
$L[f]=\partial_tf-tr_gD^2f$ and $\omega=Dw/v$, a vector field in
$D(t)$.\vspace{.2cm}

\textbf{Lemma 11.1.} For the surfaces $\Sigma_t$ evolving by graph
mean curvature motion with constant contact angle, the mean
curvature satisfies:
$$\left \{
\begin{array}{l}L[H]=|h|_g^2H+Hh^2(\omega,\omega)-H^2h(\omega,\omega)\quad \mbox{ on }D(t)\\
d_nH=\frac{\beta^2}{\beta_0}Hh_{nn}\quad \mbox{ on }\partial
D(t)\end{array}\right .$$ \vspace{.3cm}

\emph{Boundary conditions for $h_{ij}$.} Fix $p\in \partial D(t)$
and let $(\tau_a)$ be an orthonormal frame for $T_p\partial D(t)$
(in the induced metric), satisfying
$\nabla^{\Gamma}_{\tau_a}\tau_b(p)=0$ ($\nabla^{\Gamma}$ is the
connection induced on $\Gamma_t$ by the euclidean connection $d$,
or, equivalently, by $\nabla$, the Levi-Civita connection of the
metric $g$ in $D(t)$); we extend the $\tau_a$ to a tubular
neighborhood of $\Gamma_t$ so that $\bar{\nabla}_n\tau_a=0$.
Differentiating $h(n,\tau_b)=0$ along $\tau_a$, we find:
$$(\nabla_{\tau_a}h)(n,\tau_b)=-h(\nabla_{\tau_a}n,\tau_b)-h(n,\nabla_{\tau_a}\tau_b).\qquad (11.3)$$
The second fundamental form ${\cal K}(\tau,\tau')$ of $\Gamma_t$ in
$(D(t),eucl)$ (equivalently, in $(D(t),g)$) is defined by:
$$d_{\tau_a}\tau_b=\nabla^{\Gamma}_{\tau_a}\tau_b+{\cal
K}(\tau_a,\tau_b)n\quad \mbox{ on }\partial D(t).$$ To relate $\cal
K$ to $h_{|\partial D(t)}$, note that since $w=0$ on $\partial
D(t)$:
$$h(\tau_a,\tau_b)=\frac 1vd^2w(\tau_a,\tau_b)=\frac
1v(\tau_a(\tau_bw)-Dw\cdot d_{\tau_a}\tau_b) =-d_{\tau_a}\tau_b\cdot
\frac{Dw}v=-\beta_0{\cal K}(\tau_a,\tau_b).$$ (So we see that
$\Gamma_t$ convex with respect to $n$ corresponds to $\Sigma_t$
\emph{concave} over $D(t)$, as expected). In the appendix we observe
that $\nabla_{\partial_i}\partial_j=(h_{ij}/v)Dw$. Then:
$$\nabla_{\tau_a}\tau_b=\tau_a^i((\tau_b^j)_i\partial_j+\tau_b^j\nabla_{\partial_i}{\partial_j})
=d_{\tau_a}\tau_b+\frac 1v\tau_a^i\tau_b^jh_{ij}Dw$$
$$=\nabla^{\Gamma}_{\tau_a}\tau_b+{\cal
K}(\tau_a,\tau_b)n+\frac{w_n}vh(\tau_a,\tau_b)n=(-\frac
1{\beta_0}+\beta_0)h(\tau_a,\tau_b)n=-\frac{\beta^2}{\beta_0}h(\tau_a,\tau_b)n$$
at $p$, given our assumption $\nabla^{\Gamma}_{\tau_a}\tau_b(p)=0$.
We use this immediately to compute, at $p$:
$$\nabla_{\tau_a}n=\langle
\nabla_{\tau_a}n,\tau_b\rangle_g\tau_b=-\langle
n,\nabla_{\tau_a}\tau_b\rangle_g\tau_b=\frac{\beta^2}{\beta_0}|n|^2_gh(\tau_a,\tau_b)\tau_b=\frac
1{\beta_0}h(\tau_a,\tau_b)\tau_b,$$ since
$|n|^2_g=g_{ij}n^in^j=1+w_n^2=\beta^{-2}$ at $p$. Using these
expressions for $\nabla_{\tau_a}n$ and $\nabla_{\tau_a}\tau_b$ in
(11.3) above and recalling the Codazzi equations:
$$(\nabla_nh)(\tau_a,\tau_b)=(\nabla_{\tau_a}h)(n,\tau_b)=-\frac
1{\beta_0}\sum_ch(\tau_a,\tau_c)h(\tau_c,\tau_b)+\frac{\beta^2}{\beta_0}h(\tau_a,\tau_b)h_{nn}.$$
This can also be written in the form:
$$\beta_0(\nabla_nh)(\tau,\tau')=-(h^{tan})^2(\tau,\tau')+\beta^2h_{nn}h(\tau,\tau').\quad (11.4)$$
It turns out the expression for the $n$-directional derivative of
$h(\tau,\tau')$ is exactly the same (at $\partial D(t)$):
$$\beta_0d_n(h(\tau,\tau'))=-(h^{tan})^2(\tau,\tau')+\beta^2
h_{nn}h(\tau,\tau').\quad (11.5)$$ The reason is that
$\nabla_n\tau_a=0$ at the boundary, also for the $g$-connection:
$$\nabla_n\tau_a=d_n(\tau_a)+n^i\tau_a^j\nabla_{\partial_i}\partial_j=0+\frac
1vh(n,\tau_a)Dw=0,$$ so in fact:
$$(\nabla_nh)(\tau_a,\tau_b)=n(h(\tau_a,\tau_b))=d_n(h(\tau_a,\tau_b)).$$

  As done in \cite{Stahl}, we combine this with the
result for $H_n$ to compute $(\nabla_nh)(n,n)$. From:
$$H_n=\nabla_n(tr_gh)=tr_g(\nabla_nh)=\beta^2(\nabla_nh)(n,n)+\sum_a(\nabla_nh)(\tau_a,\tau_a).$$
Here we used $|n|^2_g=\beta^{-2}$ on $\partial D(t)$, which also
implies $H=\beta^2h_{nn}+\sum_ah(\tau_a,\tau_a)$. Using also
$|h^{tan}|^2=\sum(h^{tan})^2(\tau_a,\tau_a)$,  we find for
$(\nabla_nh)(n,n)$:
$$\beta^2(\nabla_nh)(n,n)=\frac{\beta^2}{\beta_0}Hh_{nn}+\frac
1{\beta_0}|h^{\tan}|^2-\frac{\beta^2}{\beta_0}(H-\beta^2h_{nn})h_{nn}$$
$$=\frac 1{\beta_0}(|h^{tan}|^2+\beta^4h_{nn}^2)=\frac
1{\beta_0}|h|^2_g,$$ since $g^{nn}=\beta^2$ at $\partial D(t)$.
Equivalently:
$$\beta_0(\nabla_nh)(n,n)=\frac 1{\beta^2}|h|^2_g\quad \mbox{ on
}\partial D(t).$$\vspace{.2cm}

It is easy to obtain the corresponding expression for the euclidean
connection. Noting that at $\partial D(t)$:
$$\nabla_nn=d_nn+n^in^j\frac 1vh_{ij}Dw=\beta_0h_{nn}n,$$
we find:
$$(d_nh)(n,n)=n(h_{nn})=(\nabla_nh)(n,n)+2h(\nabla_nn,n)=(\nabla_nh)(n,n)+2\beta_0h_{nn}^2,$$
so that:
$$\beta_0d_n(h(n,n))=\frac 1{\beta^2}|h|^2_g+2\beta_0^2h_{nn}^2\quad \mbox{ on
}\partial D(t).$$\vspace{.2cm}

We record these results as a lemma, including also the evolution
equations derived in the appendix.\vspace{.2cm}

\textbf{Lemma 11.2.} Under graph mean curvature motion with constant
contact angle, the second fundamental form satisfies the following
tensorial evolution equations. ($C_{ij}$ and $\bar{C}_{ij}$ are
symmetric 2-tensors cubic in $h$, given explicitly in the appendix.)

(i) For the operator $L=L_g$:

$$L[h_{ij}]=-2[h_i^kd_{\omega}(h_{jk})+h_j^kd_{\omega}(h_{ik})]+\bar{C}_{ij}\mbox{
on }D(t),$$ with boundary conditions on $\partial D(t)$:
$$\left \{ \begin{array}{ll}
h(n,\tau)=0\\
\beta_0d_n(h(\tau,\tau'))=-(h^{tan})^2(\tau,\tau')+\beta^2
h_{nn}h(\tau,\tau')\\
\beta_0d_n(h(n,n))=\frac 1{\beta^2}|h|^2_g+2\beta_0^2h_{nn}^2
\end{array}\right .$$

(ii) For the operator $\partial_t-\Delta_g$, where $\Delta_g$ is the
Laplace-Beltrami operator of $g$:

$$(\partial_t-\Delta_g)[h]_{ij}=H(\nabla_{\omega}h)_{ij}
+H_ih(\omega,\partial_j)+H_jh(\omega,\partial_i)+{C}_{ij}\mbox{ on
}D(t),$$ with boundary conditions on $\partial D(t)$:
$$\left \{ \begin{array}{ll}
h(n,\tau)=0\\
\beta_0(\nabla_nh)(\tau,\tau')=-(h^{tan})^2(\tau,\tau')+\beta^2h_{nn}h(\tau,\tau')\\
\beta_0(\nabla_nh)(n,n)=\frac 1{\beta^2}|h|^2_g
\end{array}\right.$$
\vspace{.2cm}

It is also useful to compute the boundary condition for $|h|^2_g$.
Using Lemma 11.2(ii), we have at $\partial D(t)$:
$$\begin{array}{l}(\beta_0/2)d_n|h|^2_g=
\beta_0\langle \nabla_nh,h\rangle_g\\
=\beta_0\beta^4(\nabla_nh)(n,n)h_{nn}+\beta_0\sum_{b,c}(\nabla_nh)(\tau_a,\tau_b)h(\tau_a,\tau_b)\\
={\beta^2}|h|^2_gh_{nn}+\sum_{a,b}[-(h^{tan})^2(\tau_a,\tau_b)+\beta^2h_{nn}h(\tau_a,\tau_b)]h(\tau_a,\tau_b)\\
=\beta^2(|h|_g^2+|h^{tan}|^2_g)h_{nn}-tr_g(h^{tan})^3\end{array}$$
Since on $\partial D(t)$:
$tr_gh^3=\beta^6(h_{nn})^3+tr_g(h^{tan})^3$, we may state this in a
slightly different form. Including also the evolution equation for
$|h|^2_g$ (see appendix), we have the following lemma.\vspace{.2cm}

\textbf{Lemma 11.3.} Under graph mean curvature flow, the function
$|h|^2_g$ satisfies the evolution equation and Neumann boundary
condition:
$$\left \{ \begin{array}{l}
(\partial_t-\Delta_g)|h|^2_g=-2|\nabla
h|_g^2+Hd_{\omega}|h|^2_g+2|h|^4_g-4Hh^3(\omega,\omega)-2H|h|^2_gh(\omega,\omega)\\
(\beta_0/2) d_n|h|^2_g=2\beta^2|h|_g^2h_{nn}-tr_g(h^3)\mbox{ on
}\partial D(t)
\end{array} \right .$$\vspace{.5cm}

\textbf{12. A maximum principle for symmetric
2-tensors.}\vspace{.2cm}

By the local existence theorem, for suitable initial data we have a
mean curvature motion $F=[\varphi,u]\in
C^{2+\alpha,1+\alpha/2}(Q_0,\mathbb{R}^{n+1})$, where $Q_0=D_0\times
[0,T]$ and, for each $t\in [0,T]$, $\varphi_t:D_0\rightarrow D(t)$
is a $C^{2+\alpha}$ diffeomorphism. In particular, with
$\delta=\alpha^2$, $w_t=u_t\circ \varphi_t^{-1}:D(t)\rightarrow
\mathbb{R}$ defines a graph m.c.m. $w\in
C^{2+\delta,1+\delta/2}(E;\mathbb{R})$ in an open spacetime domain
$$E=\bigcup_{t\in (0,T)}D(t)\times \{t\}\subset \mathbb{R}^n\times
\mathbb{R}.$$ We have a $C^{2+\alpha,1+\alpha/2}$ diffeomorphism:
$$\Phi:\bar{Q_0}\rightarrow \bar{E},\quad
\Phi(x,t)=(\varphi_t(x),t),$$ which, for any $t_0>0$, restricts to a
diffeomorphism $Q_{t_0}\rightarrow E_{t_0}$, where:
$$Q_{t_0}=D_0\times (t_0,T),\qquad E_{t_0}=\bigcup_{t\in
(t_0,T)}D(t)\times \{t\}.$$ The parabolic boundary of $E$ is the
disjoint union of `base' and `lateral boundary':
$$\partial_pE=(\bar{D}_0\times \{0\})\sqcup
\partial_lE,\quad \partial_lE=\bigcup_{t\in (0,T)}\partial D(t)\times \{t\}.$$
(The notions of `parabolic boundary, `base' and `lateral boundary'
have general definitions for arbitrary bounded spacetime domains
(see \cite{Lieberman}), but using $\Phi$ it is easy to see that they
are given by the above sets.) In particular, note that $\Phi$
defines a diffeomorphism ($C^{k+\alpha,(k+\alpha)/2)}$ up to the
lateral boundary, if $D_0$ is a $C^{k+\alpha}$ domain and $F\in
C^{k+\alpha,(k+\alpha)/2)}(Q_0)$):
$$Q_{t_0}\cup \partial_lQ_{t_0}\rightarrow E_{t_0}\cup
\partial_lE_{t_0},$$
for each $t_0>0$.\vspace{.2cm}

Denote by $L$ the operator
$L=\partial_t-g^{ij}(Dw)\partial_i\partial_j$, so $Lw=0$ in $E$ and
$w=0$ on $\partial_lE$. The following height bound is immediate.

\textbf{Lemma 12.1.} Assume $0<w_0<M$ in $D_0$.  Then $0\leq w\leq
M$ in $\bar{E}$ (and vanishes only on $\partial_lE$).

\emph{Proof.} Follows from the weak maximum principle for the
operator $L$, since $0\leq w\leq M$ holds on the parabolic boundary
$\partial_pE$. \vspace{.2cm}

It is well-known that the function $v=\sqrt{1+|Dw|^2}$ solves the
evolution equation (assuming $Dw\in C^{2,1}(\bar{E})$, see e.g.
\cite{Guan}):
$$L[v]+\frac 2vg^{ij}v_iv_j=-v|h|^2_g,\mbox{ or }L[v]=-\frac 2v|Dv|^2_g-v|h|^2_g.$$
From the maximum principle, we have the following global bound on
$v$ (equivalently, on $|Dw|$).\vspace{.2cm}

\textbf{Lemma 12.2} Assume $w$ is a solution with $Dw\in
C^{2,1}(\bar{E})$. Then we have on $\bar{E}$:
$$v(z)\leq \max\{sup_{D(t_0)}v(x,t_0),\frac 1{\beta}\}.$$

\emph{Proof.} By the weak maximum principle,
$\max_{\bar{E}}v=\max_{\partial_pE}v$. Note $v_{|S}\equiv \frac
1{\beta}.$\vspace{.2cm}

It follows from this lemma that $g_{ij}(t)$ is uniformly equivalent
to the euclidean metric in $D(t)$: if $v\leq \bar{v}$ in $\bar{E}$,
and $X$ is a vector field in $D(t)$:
$$|X|_e^2\leq |X|_g^2=g_{ij}X^iX^j=|X|^2_e+(X\cdot Dw)^2\leq
|X|^2_e(1+|Dw|^2)\leq \bar{v}^2|X|_e^2.$$ Also, if
$\omega:=v^{-1}Dw$:
$$|\omega|^2_e=\frac{|Dw|^2_e}{v^2}=1-\frac 1{v^2}\leq
1-\frac1{\bar{v}^2}.$$\vspace{.2cm}

The main result in this section is a maximum principle for symmetric
2-tensors satisfying a parabolic equation on a spacetime domain such
as $E$ (image of a cylinder under a diffeomorphism of the special
type $\Phi$). \vspace{.2cm}

Recall the `boundary point lemma' for scalar equations . It holds
for open spacetime domains $\Omega\subset \mathbb{R}^n\times
\mathbb{R}_+$ satisfying an `interior ball condition':

\emph{Interior ball condition:} For each $P=(p,\bar{t})\in
\partial_l\Omega$ we may find a ball $B$ (in the euclidean metric in
$\mathbb{R}^{n+1})$ which is tangent to $\partial_l \Omega$ only at
$P$ and satisfies: (i) the line segment from $P$ to the center of
the ball is not parallel to the $t$ axis; (ii) $B\cap \{t\leq
\bar{t}\}\subset \Omega \cap \{t\leq \bar{t}\}$. \vspace{.2cm}

\emph{Remark.} For the domain of interest the interior ball
condition follows from the fact that $\partial_lE=\Phi(\partial
D_0\times (0,T))$, with $\Phi\in C^{2,1}(D_0\times (0,T))$ of the
special form above.\vspace{.2cm}

\textbf{Lemma 12.3.}(\cite{ProtterWeinberger} Thm.6, p.174.) Let
$\Omega\subset \mathbb{R}^n\times \mathbb{R}_+$ be a connected open
set satisfying the interior ball condition. Assume $f\in
C^{2,1}(\Omega)$ satisfies the uniformly parabolic inequality:
$$\partial_tf-tr_gd^2f-d_Xf\leq 0.$$
Here $g=g_t$ is a Riemannian metric in each section $\Omega(t)$, and
$X_t$ is a bounded vector field in $\Omega(t)$. Denote by $n=n_t$
the inner unit normal of $\Omega(t)$.\vspace{.2cm}

Assume the supremum $M$ of $f$ in $\Omega_{\bar{t}}:=\Omega \cap
\{t\leq \bar{t}\}$ is attained at the point $P\in \partial
\Omega(\bar{t})$, and that $f<M$ for $t<\bar{t}$. Then
$d_nf(P)<0$.\vspace{.2cm}

We now state the hypotheses of our tensorial maximum principle.

$E\subset \mathbb{R}^n\times [0,T]$ is the image of a cylinder
$D_0\times (0,T)$ under a $C^{3,2}$ diffeomorphism $\Phi$ of the
form $\Phi(x,t)=(\varphi_t(x),t)$, with $\varphi_t:D_0\rightarrow
D(t)$ a $C^3$ diffeomorphism up to the boundary, for each $t\in
[0,T]$ ($D(t)$ is the $t=const.$ section of $E$); $\bar{D_0}\subset
\mathbb{R}^n$ is assumed to be the image of the closed unit ball
under a $C^3$ diffeomorphism.) In particular, the lateral boundary
$\partial_lE$ is of class $C^{3,2}$. On $\partial_lE$ we have the
\emph{inner} unit normal $n=n_t\in \mathbb R^n$. Extend $n_t$ to a
vector field in all of $\bar{D}(t)$ so that it is in
$C^{2,1}(\bar{E},\mathbb R^n)$, arbitrarily except for the
requirements that $|n|\leq 1$ pointwise and $d_nn=0$ in a tubular
neighborhood of $\partial D(t)$ (equivalently, $n^i\partial_in^j=0$
for each $j$.) Fix $R>0$ so that $D(t)\subset B_R(0)$, for each
$t\in [0,T]$\vspace{.2cm}

 The assumptions on the coefficients are given next.

 $g=g_t$ is a $t$-dependent Riemannian metric in $\bar{D}(t)$, uniformly equivalent to
 the euclidean metric for $t\in [0,T]$;

 $X=X_t$ is a bounded $t$-dependent vector field in $\bar{D}(t)$;

 $q=q(z,m)$ assigns to each $z\in \bar{E}$ and each $m$ in
 $\mathbb{S}$ (the space of quadratic forms in $\mathbb R^n$) a quadratic form $q\in
 \mathbb{S}$. $q$ is assumed to be $C^{2,1}$ in $z$, locally Lipschitz in
 $m$ (uniformly in $z\in \bar{E}$);

 $b=b(z,m)\in \mathbb{S}$ is defined for $z\in \partial_lE$, with the
 same regularity assumptions as $q$.\vspace{.3cm}

 \emph{Remark.} We state the following theorem in terms of the
 Laplace-Beltrami heat operator $\partial_t-\Delta_g$ and the $g$- Riemannian connection
 $\nabla$, but the result also holds for $L$ and the `euclidean
 connection' $d$.\vspace{.2cm}

 \textbf{Theorem 12.1.} Assume $m\in C^{2,1}(\bar{E};\mathbb{S})$ satisfies in
 $E$ the tensorial differential inequality:
 $$\partial_tm_{ij}-(\Delta_g m)_{ij}\leq (\nabla_Xm)_{ij}
 +q_{ij}(\cdot,m(\cdot)),$$
 and on $\partial_lE$ the boundary condition:
 $$(\nabla_nm)_{ij}(z)\geq b_{ij}(z,m(z)).$$
 Suppose the functions $q$ and $b$ satisfy the following `null
 eigenvector conditions': if, for some $\hat{m}\in \mathbb{S}$,  $V\in {\mathbb R}^n$ is a null
 eigenvector of $\hat{m}$ ($\hat{m}_{ij}V^j=0\forall i$), then, for
 any $z\in \bar{E}$ (resp. any $z\in \partial_lE$):
 $$q_{ij}(z,\hat{m})V^iV^j\leq 0\quad \mbox{ (resp. }b_{ij}(z,\hat{m})V^iV^j\geq
 0).$$
 Then weak concavity of $m$ at $t=0$ is preserved:
 $$m\leq 0\mbox{ in } D(0)\Rightarrow m\leq 0\mbox{ in }\bar{E}.$$
 \vspace{.3cm}

 \emph{Proof.} The assumptions imply there is $K>0$ (depending only on $E$ and
 on the functions $X$, $g$, $n$, $q$ and $b$) satisfying:

 $$|n|_{C^{2,1}(\bar{E})}\leq K,\quad |X(z)|_{eucl}\leq K,\quad |g(z)|+|g^{-1}(z)|\leq K,\quad z\in
 \bar{E};$$
 and if $m,\hat{m}\in C^{2,1}(\bar{E},\mathbb{S})$ satisfy (for some
 $\mu:\bar{E}\rightarrow \mathbb{R}_+$):
 $$-\mu(z)g\leq m(z)-\hat{m}(z)\leq \mu(z)g$$
 (where the inequality of quadratic
 forms has the usual meaning), then also:
 $$q(z,m(z))\leq q(z,\hat{m}(z))+K\mu(z)g,\quad z\in \bar{E},$$
 $$b(z,m(z))\geq b(z,\hat{m}(z))-K\mu(z)g,\quad z\in
 \partial_lE.$$\vspace{.2cm}

 Now define, for $z\in \bar{E}$, $z=(x,t)$:
 $$\varphi(z):=-2Kn(z)\cdot x:=2Ks(z),$$
 where we use the euclidean inner product and, on
 $\partial_lE$, $s$ is the `support function' of $\partial D(t)$
 (positive if $D(t)$ is convex and contains the origin). It is clear we may find $M=M(R,K)>0$
 depending only on $K, R$ and $|n|_{C^{2,1}}$ so that:
 $$|\varphi|_{C^{2,1}}\leq M, \quad
 |d\varphi|_g^2+|\Delta_g\varphi|\leq M,\quad |X\cdot d\varphi|\leq
 M.$$
 We assume also $M\geq K$. Now, given $m$ as in the statement of the theorem and
 given constants $\epsilon>0, \gamma>0$ and
 $\delta >0$, define for $ E^{\delta}:=E\cap \{t<\delta\}$:
 $$\hat{m}(z):=m(z)-(\epsilon t+\gamma
 e^{\varphi(z)})g,\quad
 z\in \bar{E}^\delta.$$
 Clearly $\hat{m}\in C^{2,1}(\bar{E}^{\delta}; \mathbb{S})$. We now
 derive the constraints on $\delta$, $\epsilon$ and $\gamma$. It
 will turn out that $\delta$ must be taken small enough (depending only on
 $K,R$), $\epsilon>0$ is arbitrary and $\gamma$ is $\epsilon$ times
 a constant depending only on $K,R$.\vspace{.2cm}

 The following inequalities are easily derived:
 $$q(z,m(z))\leq q(z,\hat{m}(z))+K(\epsilon t+\gamma
 e^{\varphi(z)})g;$$
 $$\nabla_Xm=\nabla_X\hat{m}+\gamma (e^{\varphi}
 d_X\varphi)g\leq \nabla_X\hat{m}+(\gamma
 e^{\varphi}M)g;$$
 $$\partial_t\hat{m}=\partial_tm-\epsilon g-(\gamma
 e^{\varphi}\partial_t\varphi)g
 \leq \partial_tm+(\gamma
 e^{\varphi}M)g-\epsilon g;$$
 $$\Delta_g\hat{m}=\Delta_gd^2\hat{m}-\gamma
 e^{\varphi}(|d\varphi|^2_g+\Delta_g\varphi)g
 \geq \Delta_gm-(\gamma e^{\varphi}M)g.$$
 $$b(z,m(z))\geq b(z,\hat{m}(z))-K(\epsilon t+\gamma
 e^{\varphi})g$$
 \vspace{.2cm}

 We use this to compute:
 $$\partial_t\hat{m}-\Delta_g\hat{m}\leq \partial_t
 m-\Delta_gm+(2\gamma e^{\varphi}M)g-\epsilon g$$
 $$\leq q(z,m(z))+\nabla_Xm+(2\gamma e^{\varphi}M)g-\epsilon
 g$$
 $$\leq q(z,\hat{m}(z))+\nabla_X\hat{m}+K(\epsilon
 t+\gamma e^{\varphi})g+(3M\gamma e^{\varphi})g-\epsilon g$$
 $$\leq q(z,\hat{m}(z))+\nabla_X\hat{m}+M\epsilon
 tg+4M\gamma e^{\varphi}g-\epsilon g$$
 (using $K\leq M$ in the last step). We conclude the inequality:
 $$\partial_t\hat{m}-\Delta_g\hat{m}\leq q(z,\hat{m}(z))+\nabla_X\hat{m}
 -({\epsilon}/2)g\qquad (13.1)$$ will hold in
$E^{\delta}$, provided the constants are selected
 so that, for $z\in E^{\delta}$:
 $$4M\gamma e^{\varphi(z)}+M\epsilon t\leq {\epsilon}/2.\qquad
 (A)$$\vspace{.2cm}
 Turning to boundary points $z=(x,t)\in \partial_lE$, note that
 $d_n\varphi=-2K$, so that:
 $$\nabla_n\hat{m}(z)=\nabla_nm(z)-(\gamma
 e^{\varphi(z)}d_n\varphi(z))g\geq b(z,m(z))-(\gamma
 e^{\varphi(z)}d_n\varphi(z))g$$
 $$\geq b(z,\hat{m}(z))-K(\epsilon t+\gamma
 e^{\varphi(z)})g-(\gamma
 e^{\varphi(z)}d_n\varphi(z))g$$
 $$\geq b(z,\hat{m}(z))+K(\gamma e^{\varphi(z)}-\epsilon
 t)g,$$
 implying the inequality:
 $$\nabla_n\hat{m}(z)\geq b(z,\hat{m}),\quad z\in \partial_lE^{\delta}\qquad (13.2)$$
 will hold provided the constants are chosen so that, on
 $\partial_lE^{\delta}$:
 $$\epsilon t\leq \gamma e^{\varphi(z)}.\qquad (B)$$
 Bearing in mind that, on $E$: $e^{-2KR}\leq e^{\varphi(z)}\leq
 e^{2KR}$, it is not hard to arrange for (A) and (B) to hold, or
 equivalently, for:
 $$\epsilon t\leq \gamma e^{\varphi(z)},\qquad 10M\gamma
 e^{\varphi(z)}\leq \epsilon.$$
 Given $\epsilon >0$, define $\gamma$ so that $10M\gamma
 e^{2KR}=\epsilon$. Then the second inequality holds, and so will
 the first, provided:
 $$\epsilon t\leq \gamma e^{-2KR}=(\epsilon/10M)e^{-4KR},$$
 which is true for any $\epsilon>0$, if $\delta$ is defined via
 $\delta:=e^{-4KR}/10M$ (recall $t\in [0,\delta]$).\vspace{.2cm}

  Note that, since $m\leq 0$ at $t=0$, it follows that $\hat{m}$ is negative-definite
  at $t=0$, and hence also
 for small time, and we \emph{claim} that this persists
 throughout $\bar{E}^{\delta}$, so that (letting $\epsilon
 \rightarrow 0$) $m\leq 0$ in $\bar{E}^{\delta}$. Restarting the
 argument at $t=\delta$, we see this is enough to prove the theorem.
 \vspace{.2cm}

 To prove this claim, suppose (by contradiction) $\hat{m}$ acquires
 a null eigenvector $0\neq V\in {\mathbb R}^n$ at a point
 $z_1=(x_1,t_1)\in \bar{E}^{\delta}$, with $t_1\in (0,\delta]$ the
 first time this happens.

 Let $\hat{f}(z):=\hat{m}_{ij}V^iV^j, \quad z\in E^{\delta}$ (that
 is, we `extend' $V$ to $E^{\delta}$ as a constant vector.) It
 follows from (13.1) that $\hat{f}$ satisfies in
 $E^{\delta}$:
 $$\partial_t\hat{f}\leq (\Delta_g\hat{m})_{ij}V^iV^j+
 (\nabla_X\hat{m})_{ij}V^iV^j+q_{ij}(\cdot,\hat{m})V^iV^j-\frac{\epsilon}2|V|^2_{g}.$$
 A short and standard Riemannian calculation (using the fact that $V$
 is a null eigenvector for $\hat{m}$) shows that:
 $$d_X\hat{f}=(\nabla_X\hat{m})_{ij}V^iV^j,\quad
 \Delta_g\hat{f}=(\Delta_g\hat{m})_{ij}V^iV^j.$$
  Using the null eigenvector
 condition for $q$, we find that $\hat{f}$ satisfies in $E^{\delta}$
 the strict inequality:
 $$\partial_t\hat{f}<tr_gd^2\hat{f}+d_X\hat{f}.$$
 This shows $x_1$ cannot be an interior point of $D(t_1)$, for then
 (as a first-time interior maximum point for $\hat{f}$) we would
 have $\Delta_g\hat{f}(z_1)\leq 0$ and $d\hat{f}(z_1)=0$,
 contradicting $\partial_t\hat{f}(z_1)\geq 0$. Thus $x_1\in \partial
 D(t_1)$. Since $\hat{f}$ satisfies the differential inequality just
 stated and $z_1=(x_1,t_1)$ is a first-time boundary maximum in
 $\bar{E}^{\delta}$, the parabolic Hopf lemma (Lemma 12.3 above) implies
 $d_n\hat{f}(z_1)<0$. On the other hand, as seen above in (13.2):
 $$d_n\hat{f}= (\nabla_n\hat{m})_{ij}V^iV^j\geq
 b_{ij}(z_1,\hat{m}(z_1))V^iV^j\geq 0,$$
 from the boundary null-eigenvector condition. This contradiction
 concludes the proof.\vspace{.3cm}

 \textbf{Corollary 12.2.} Suppose $m\in C^{2,1}(\bar{E},\mathbb{S})$
 satisfies the same differential inequality, with the same
 hypotheses on the coefficients as in the theorem (including the
 null eigenvector condition for $q$), and the boundary conditions:
 $$\left \{ \begin{array}{l}
 m(z)(n,\tau)=0,\quad \forall z=(x,t)\in \partial_lE,\tau\in
 T_x\partial D(t)\\
 ({\nabla}_nm)(n,n)\geq b_{nn}(z,m(z))\\
 ({\nabla}_nm)(\tau,\tau)\geq
 b^{tan}(z,m(z))(\tau,\tau),\quad \tau \in T_x\partial D(t),\\
 \end{array}\right .$$
 for functions $b_{nn}(z,\hat{m})$ from $\partial_lE\times \mathbb{S}$ to $\mathbb{R}$ and $b^{tan}$ assigning
 to $(z,\hat{m})\in \partial_lE\times \mathbb{S}$, $z=(x,t)$, a quadratic form in $T_{x}\partial D(t)$.
 Suppose $b_{nn}\geq 0$ in $E\times \mathbb{S}$ and $b^{tan}$
 satisfies for each
 $\hat{m}\in \mathbb{S}$:
 $$ \hat{m}_{ij}\tau^i=0\mbox{ for some }\tau\in T_x\partial
 D(t)\Rightarrow b^{tan}(z,\hat{m})(\tau,\tau)\geq 0.$$
 Then, as in the theorem, weak concavity is preserved:
 $$m\leq 0\mbox{ at }t=0\Rightarrow m\leq 0\mbox{ in
 }\bar{E}.$$\vspace{.2cm}

 \emph{Proof.} This is proved as the theorem, with the following
 change in the last part of the proof: if $0\neq V\in \mathbb{R}^n$
 is a null eigenvector of $\hat{m}$ (defined as in the proof of the
 theorem) at a boundary point $z_1=(x_1,t_1)\in \partial_l E$,
 write:
 $$V=V^nn+V^T,\quad V^T\in T_{x_1}\partial D(t_1).$$
 Assume first $V^n\neq 0$. Then (noting that $\hat{m}$ splits at the
 boundary if $m$ does), we see that $n$ is a null eigenvector of
 $\hat{m}$ at $z_1$, so we define
 $\hat{f}(z)=\hat{m}_{ij}(z)n^i(z_1)n^j(z_1)$ and repeat the
 argument. At $z_1$,
 $({\nabla}_n\hat{m})(n,n)\geq b_{nn}(z_1,\hat{m}(z_1))\geq 0$ leads to a
 contradiction with the parabolic Hopf lemma, as before.

 If $V^n=0$, then $V^T\in T_{x_1}\partial D(t_1)$ must be a null eigenvector of $\hat{m}$ at
 the boundary point $z_1$, and then we run the argument with
 $\hat{f}(z)=\hat{m}(z)(V^T,V^T)$, leading to a contradiction,
 as before.\vspace{.3cm}

 \textbf{Corollary 12.3.} Let $w\in C^{4,2}(E)$ define a MCM of graphs with constant-angle
 boundary conditions, where $E$ is as in the statement of theorem 12.1. Then weak concavity is preserved:
 $$h\leq 0\mbox{ at }t=0\Rightarrow h\leq 0\mbox{ in }\bar{E}.$$

 \emph{Proof.}  From Lemma 11.2, $h$ satisfies:
 $$(\partial_t-\Delta_g)h_{ij}=H\nabla_{\omega}h_{ij}+q(z,h)_{ij},$$
 $$q(z,h)_{ij}=H_ih(\omega,\partial_j)+H_jh(\omega,\partial_i)+|h|_g^2h_{ij}
 +Hh(\partial_i,\omega)h(\partial_j,\omega)-Hh(\omega,\omega)h_{ij}.$$
 (Here $H_i,H_j,H$ and $\omega^i$ are regarded as fixed functions of
 $z\in E$.) Clearly $q$ satisfies the null eigenvector condition,
 since $q_{ij}V^iV^j=0$ when $h_{ij}V^j=0$ for all $i$. In addition,
 expressions obtained for $d_nh$ in lemma 11.2 show that the
 boundary conditions in Corollary 12.2 are satisfied, with:
 $$b_{nn}(z,\hat{m})\equiv 0,\quad
 b^{tan}(z,\hat{m})=-((\hat{m})^{tan})^2+\beta^2\hat{m}_{nn}\hat{m}^{tan}.$$
  Hence the
 claim follows from Corollary 12.2.\vspace{.3cm}

 \emph{Remark 12.1.} For less regular solutions, we may apply the
 theorem to a domain $E_{t_0}=E\cap \{ t>t_0\}$, for arbitrarily
 small $\delta>0$. Thus, assuming $h<0$ at $t=0$ (strictly
 negative-definite), we conclude from corollary 13.3 that $h\leq 0$
 for all $t$.\vspace{.2cm}

 \emph{Remark 12.2.} It seems plausible that a slightly different version
 of the result in this section could be used to strengthen the
 conclusions of \cite{Stahl}.\vspace{.2cm}

 \emph{Finite existence time.}\vspace{.2cm}

It is not difficult to derive that the flow is defined only for
finite time in the concave case.\vspace{.2cm}

\textbf{Lemma 12.4.} Let $w\in C^{4,2}(E), E\subset {\mathbb
R}^n\times [0,T)$ define a graph MCM $\Sigma_t$ with constant-angle
boundary conditions on a moving boundary. Assume $\Sigma_0$ (and
hence $\Sigma_t$, for all $t$) is weakly concave. Then:

Assume $H_{|t=0}\leq H_0<0$ (where $H_0$ is a negative constant).
Then $T\leq t_*=\frac 1{2H_0^2c_n}$ ( we are assuming $T=\sup \{t\in
[0,T); D(t)\neq \emptyset\}$). Here $c_n>0$ depends only on $n$ and
an upper bound for $v$ in $E$. \vspace{.2cm}

The proof is based on the evolution equation and boundary condition
for $H$ (see Appendix 2: $\omega=Dw/v$):
$$L[H]=|h|^2_gH+Hh^2(\omega,\omega)-H^2h(\omega,\omega),\quad
H_n=(\beta^2/\beta_0)Hh_{nn}.$$
 Since $h^2(\omega,\omega)\geq 0$,  $|h|^2_g \geq
(1/n)H^2$ and (given that $h\leq 0$) $h(\omega,\omega) \geq
|Dw|^2H$, we have:
$$L[H]\leq \frac 1nH^3+|Dw|^2H^3\leq c_nH^3,$$
where $c_n$ depends on $n$ and on $\sup_E |v|$ (already known to be
finite).  Let $\phi(t)$ solve the o.d.e.
$\dot{\phi}=c_n\phi^3,\phi(0)=H_0$:
$$\phi(t)=H_0[1-2c_nH_0^2t]^{-1/2},\quad 0\leq t<t_*:=\frac 1{2H_0^2c_n}.$$
Then with $\psi:=(1/n)(H^2+H\phi+\phi^2)>0$, setting $\chi=H-\phi$:
$$L[\chi]\leq\psi \chi\quad \mbox{ in }E;$$
$$\chi_n=\frac{\beta^2}{\beta_0}(\chi+\phi)h_{nn}\geq
\frac{\beta^2}{\beta_0}\chi\quad \mbox{ on }\partial_lE$$ (since
$\phi<0$ and $h_{nn}\leq 0$). Given that $\chi\leq 0$ at $t=0$, it
follows from the maximum principle that $\chi\leq 0$, or $H\leq
\phi$ in $[0,\min\{T,t_*\})$. This shows $t_*<T$ is impossible,
since $\phi \rightarrow -\infty$ as $t\rightarrow t_*$.
\vspace{.2cm}

\emph{Remark 12.1.} It would be natural to try to show that a
negative upper bound $H_0$ on the mean curvature (at $t=0$) is
preserved, at least under the assumption of concavity.
Unfortunately, the evolution equation for $H$ (under graph m.c.m.)
does not lend itself to a maximum principle argument. Letting
$u:=H-H_0$, we have

$$L[u]= |h|^2_gu+uh^2(\omega,\omega)-u(H+H_0)h(\omega,\omega)+H_0Q\quad \mbox{ in
}E,$$
$$Q:=|h|_g^2+h^2(\omega,\omega)-H_0h(\omega,\omega).$$ At a
point where $u=0$, we would need to show $L[u]\leq 0$. But it is not
true that $Q\geq 0$ at such a point. ($u_n\geq 0$ does hold at
boundary points.)\vspace{.2cm}

The exception is if $n=2$ (under an additional condition). Let
$\hat{\omega}=\omega/|\omega|_g,
\tilde{\omega}={\omega}^{\perp}/|\omega|_g$. It is easy to check
that ${\cal B}=\{\hat{\omega},\tilde{\omega}\}$ is a $g$-orthonormal
frame at each point where $\omega\neq 0$. Then with:
$$a:=h(\hat{\omega},\hat{\omega}),\quad
b:=h(\hat{\omega},\tilde{\omega}),\quad c:=h(\tilde{\omega},
\tilde{\omega}),$$ we have:
$$h^2(\hat{\omega},\hat{\omega})-Hh(\hat{\omega},\hat{\omega})=a^2+b^2-(a+c)a=b^2-ac=-\Delta,$$
where $\Delta$, the determinant of the matrix of $h$ in $\cal B$, is
non-negative if $h\leq 0$. In particular:
$$h^2(\omega,\omega)-Hh(\omega,\omega)=-|\omega|^2_g\Delta\leq 0$$
in the concave case. Now consider the expression denoted above by
$Q$ (at a point where $u=0$, or $H=H_0$). Since
$|\omega|_g^2=|Dw|^2$:
$$Q=|h|_g^2+h^2(\omega,\omega)-Hh(\omega,\omega)=a^2+2b^2+c^2+|Dw|^2(b^2-ac)$$
$$=b^2(2+|Dw|^2)+a^2-|Dw|^2ac+c^2,$$
so $Q\geq 0$ provided $|Dw|^2\leq 2$. This last condition is
equivalent to $v\leq \sqrt{3}$, and hence (Lemma 12.2) is preserved
by the evolution if it holds at $t=0$. Thus we have:\vspace{.2cm}

\textbf{Proposition 12.5.} Assume $n=2$, $h\leq 0$ and $v\leq
\sqrt{3}$ on $\Sigma_0$ (in particular, $\beta\geq 1/\sqrt{3}$).
Then $H\leq H_0 <0$ at $t=0$ implies $H\leq H_0$ for all $t\in
[0,T_{max})$.\vspace{.3cm}

 \textbf{13. Global bounds from boundary bounds for $\nabla^nh$.} \vspace{.2cm}

 In this section we begin to develop a continuation criterion for solutions of graph
 mean curvature motion with constant contact angle, based on the second
 fundamental form. Our first observation is that the
supremum of $|h|_g$ on the moving boundary controls its value in the
interior. Recall we already have a bound on $\sup_E v$ (Lemma 12.2);
it is a well known-fact for mean curvature flow of graphs that this
implies interior bounds for the second fundamental form and its
covariant derivatives (\cite{EckerHuisken}, \cite{EckerBook}). In
the next lemma we describe a global bound for mean curvature motion
of graphs with moving boundaries.\vspace{.2cm}

\textbf{Lemma 13.1.} Let $w:E\rightarrow \mathbb{R}$ be a
(sufficiently regular) solution of graph m.c.m in a spacetime domain
$E\subset \mathbb{R}^n\times [0,T]$, where $T<\infty$. Assume the
first derivative bound $v(x,t)\leq \bar{v}$ holds globally in
$\bar{E}$. Then if the bound $|h|_g\leq h_0$ holds on the parabolic
boundary $\partial_pE$, we also have the global bound:
$$|h|_g\leq a_0\quad \mbox{ in }\bar{E},$$
for a constant $a_0$ depending only on $n,\bar{v},h_0,T$ and the
initial data of $w$. \vspace{.2cm}

\emph{Proof in the concave case.} The proof is simpler under the
assumption that $h$ is negative-definite. (As shown in the previous
section, this is preserved if it holds at $t=0$). \vspace{.2cm}

\emph{Notation:} In this proof, the norms of tensors in $D(t)$ are
always taken with respect to the induced metric $g$, so we write
e.g. $|h|$ for $|h|_g$, $|\nabla h|$ for $|\nabla h|_g$, and
$|Df|^2=g^{ij}f_if_j$ for a function $f$.\vspace{.2cm}

Recall the evolution equations:
$$L[v]=-v|h|^2-\frac 2v|Dv|^2, \quad \mbox{ so
}L[v^2]=-2v^2|h|^2-6|Dv|^2,$$
$$L[|h|^2]=-2|\nabla
h|^2+2|h|^4-4Hh^3(\omega,\omega)-2H|h|^2h(\omega,\omega).$$ In the
concave case, since $H\leq 0$ and $h^3$ is negative-definite, this
implies:
$$L[|h|^2]\leq -2|\nabla h|^2+2|h|^4.$$
The idea of the proof is to apply the maximum principle to
$f=|h|^2v^2$. In the evolution equation for $f$:
$$L[f]=v^2L[|h|^2]+|h|^2L[v^2]-2\langle D|h|^2,Dv^2\rangle_g,$$
the terms $\pm 2v^2|h|^4$ cancel exactly, and we have the
inequality:
$$L[f]\leq -2v^2|\nabla h|^2-6|h|^2|Dv|^2-2\langle
D|h|^2,Dv^2\rangle_g.$$ The term with the inner product can be
estimated in two ways: 
$$|\langle D|h|^2,Dv^2\rangle_g|\leq |D|h|^2||Dv^2|\leq 4|h|v|\nabla
h||Dv|\leq 2v^2|\nabla h|^2+2|h|^2|Dv|^2$$ and:
$$\langle D|h|^2,Dv^2\rangle_g=\frac 1{v^2}\langle
D(|h|^2v^2),Dv^2\rangle_g-\frac{|h|^2}{v^2}|Dv^2|^2=\frac
1{v^2}\langle Df,Dv^2\rangle_g-4|h|^2|Dv|^2.$$ Using the second
expression, we have:
$$L[f]\leq -2v^2|\nabla h|^2-6|h|^2|Dv|^2-\frac 1{v^2}\langle
Df,Dv^2\rangle_g+4|h|^2|Dv|^2-\langle D|h|^2,Dv^2\rangle_g,$$ and
then estimating the remaining inner product term from the first
expression:
$$L[f]\leq -2v^2|\nabla h|^2-6|h|^2|Dv|^2-\frac 1{v^2}\langle
Df,Dv^2\rangle_g+4|h|^2|Dv|^2+2v^2|\nabla h|^2+2|h|^2|Dv|^2,$$
yielding after cancelation:
$$L[f]\leq -\frac 1{v^2}\langle Df,Dv^2\rangle_g.$$
Applying the (weak) maximum principle to $f$, we conclude:
$$\max_{\bar{E}}|h|^2\leq \max_{\bar{E}}f\leq \max_{\partial_pE}f\leq
\bar{v}^2\max_{\partial_pE}|h|^2,$$ which implies the result, with
an explicit constant $a_0=\bar{v}h_0.$\vspace{.2cm}

\emph{Remark.} In the general (not necessarily concave) case, we
have:
$$L[|h|^2]\leq -2|\nabla h|^2+c_n|h|^4.$$
Then the proof follows the same lines as that of proposition 3.21 in
\cite{EckerBook}: we apply the maximum principle to
$f=|h|^2(\eta\circ v^2)$, for a carefully chosen function
$\eta(s)$.\vspace{.2cm}

\emph{Evolution of $|\nabla h|^2$.}\vspace{.2cm}

In the calculation that follows, we adopt the usual convention that
in symbols such as: $$\nabla^2h*(\nabla h)^{(2)}*h^{(3)},\quad
(\nabla^jh)^{(p)}=\nabla^jh*\ldots*\nabla^jh\mbox{ ($p$ times ) },$$
$*$ denotes some unspecified $g$-contraction of the tensors in
question.\vspace{.2cm}

For the time derivative, we have:
$$\partial_t|\nabla h|^2=2\langle \partial_t(\nabla h),\nabla
h\rangle+\partial_t(g^{ij}g^{pq}g^{rs})(\nabla_ih)_{pr}(\nabla_jh)_{qs}$$
$$=2\langle \partial_t(\nabla h),\nabla
h\rangle+3(\partial_tg^{ij})\langle \nabla_ih,\nabla_jh\rangle,$$
using the Codazzi identity.

For the Hessian (using $\nabla_k\partial_l=h_{kl}\omega$, see
appendix):
$$\begin{array}{ll}\nabla^2_{k,l}|\nabla h|^2&=2\langle \nabla_l(\nabla_k\nabla
h),\nabla h\rangle+2\langle \nabla_k\nabla h,\nabla_l\nabla h\rangle
-h_{kl}d_{\omega}|\nabla h|^2\\
&=2\langle \nabla^2_{k,l}(\nabla h),\nabla h\rangle +2\langle
h_{kl}\nabla_{\omega}\nabla h,\nabla h\rangle+2\langle
\nabla_k\nabla h,\nabla_l\nabla h\rangle -h_{kl}d_{\omega}|\nabla
h|^2\\
&=2\langle\nabla^2_{k,l}(\nabla h),\nabla h\rangle+2\langle
\nabla_k\nabla h,\nabla_l\nabla h\rangle,\end{array}$$ after
cancelation. Taking traces we find:

$$(\partial_t-\Delta)|\nabla h|^2=-2|\nabla^2h|^2+2\langle (\partial_t-\Delta)(\nabla h),\nabla h\rangle
+3(\partial_tg^{ij})\langle \nabla_ih,\nabla_jh\rangle.$$

Commutation of covariant derivatives introduces the Riemann
curvature tensor, and the time derivative of the connection is also
needed:
$$(\partial_t-\Delta)(\nabla h)=\nabla
[(\partial_t-\Delta)h]+(\nabla Rm)*h+Rm*(\nabla
h)+(\partial_t\Gamma)*h,$$ where (see appendix):
$$\partial_th=\nabla dH+H\nabla_{\omega}h+T+h^{(3)},\quad
T_{ij}=H_ih(\omega,\partial_j)+H_jh(\omega,\partial_i),$$ which
combined with $\Gamma=h\omega$ and $\partial_t\omega=\nabla
H+h^{(2)}$ is easily seen to imply:
$$\partial_t\Gamma=
(\nabla dH)\omega+\nabla h*h+h^{(3)}\sim \nabla^2h+\nabla
h*h+h^{(3)}.$$ From the Gauss equation, $Rm\sim h*h$. Thus:
$$\langle(\partial_t-\Delta)(\nabla h),\nabla h\rangle\sim\langle \nabla
[(\partial_t-\Delta)h],\nabla h\rangle+\nabla^2h*\nabla h*h+(\nabla
h)^{(2)}*h^{(2)}+\nabla h*h^{(4)}.$$ On the other hand, from the
evolution equation for $h$ (appendix) we have:
$$\langle \nabla
[(\partial_t-\Delta)h],\nabla h\rangle=\langle
\nabla(H\nabla_{\omega}h+T+h^{(3)}),\nabla h\rangle= \langle \nabla
(H\nabla_{\omega}h),\nabla h\rangle +\langle \nabla T,\nabla
h\rangle +(\nabla h)^{(2)}*h^{(2)}.$$ Computing the terms on the
right, we find:
$$\langle \nabla (H\nabla_{\omega}h),\nabla h\rangle=\langle
\nabla_{\omega}h,\nabla_{\nabla H}h\rangle+H\langle
\nabla(\nabla_{\omega}h),\nabla h\rangle=\langle
\nabla_{\omega}h,\nabla_{\nabla H}h\rangle+\nabla^2h*\nabla h*h,$$
and using the Codazzi identity:
$$\langle \nabla T,\nabla h\rangle=2\langle
\nabla_{\omega}h,\nabla_{\nabla H}h\rangle +\nabla^2h*\nabla
h*h+(\nabla h)^{(2)}*h^{(2)}.$$ Putting together these results, we
have:
$$\langle (\partial_t-\Delta)(\nabla h),\nabla h\rangle=3\langle
\nabla_{\omega}h,\nabla_{\nabla H}h\rangle +\nabla^2h*\nabla
h*h+(\nabla h)^{(2)}*h^{(2)}+\nabla h*h^{(4)}.$$ On the other hand,
using the expression for $\partial_tg^{ij}$ given in the appendix we
find:
$$3\partial_tg^{ij}\langle \nabla_ih,\nabla_jh\rangle=-6\langle
\nabla_{\omega}h,\nabla_{\nabla H}h\rangle+(\nabla
h)^{(2)}*h^{(2)}.$$ So we have a cancelation, and obtain the
evolution equation:
$$(\partial_t-\Delta)|\nabla h|^2=-2|\nabla^2h|^2+
\nabla^2h*\nabla h*h+(\nabla h)^{(2)}*h^{(2)}+\nabla h*h^{(4)}.$$
\emph{Remark:} without the cancelation, the right-hand side would
involve terms of type $(\nabla h)^{(3)}$, which would be a problem
for the argument that follows.\vspace{.2cm}

Given this calculation, the following lemma has a very simple
proof.\vspace{.2cm}

\textbf{Lemma 13.2.} For a solution $w\in C^{5,3}(E)$, assume we
have a uniform bound for $h$: $|h|\leq a_0$ in $E$. Then we may find
constants $\alpha>0,C>0$ depending only on the dimension and $a_0$,
so that the function:
$$f(x,t)=\alpha |\nabla h|^2+|h|^2$$
is a subsolution in $E$: $(\partial_t-\Delta)f\leq C$.\vspace{.2cm}

\emph{Proof.} The calculation above implies:
$$(\partial_t-\Delta)|\nabla
h|^2\leq-2|\nabla^2h|^2+c_n(a_0|\nabla^2h||\nabla h|+a_0^2|\nabla
h|^2+a_0^4|\nabla h|),$$ while the evolution equation for $|h|^2$
implies:
$$(\partial_t-\Delta h)|h|^2\leq -2|\nabla h|^2+c_n (a_0^2|\nabla
h|+a_0^4).$$ Clearly we may choose $\alpha$ small enough to satisfy
the claim.\vspace{.2cm}

Our next goal is to extend this argument to higher covariant
derivatives of $h$. It turns out this does not involve a cancelation
similar to the one noted above. The various terms $T$ appearing in
each expression below are all of the same `weight', where the weight
$w$ of a term is the positive integer defined by:
$$w[\nabla^jh]=j+1,\quad w[(\nabla^jh)^{(p)}]=p(j+1),\quad j\geq 0,
p\geq 1,$$
$$\mbox{ If }T=(\nabla^{j_1}h)^{(p_1)}*\ldots*(\nabla^{j_r}h)^{(p_r)},\quad w[T]=\sum_{i=1}^rp_i(j_i+1).$$
We introduce a convenient notation for the `error terms'. For
integers $w_0\geq 1$ and $n\geq $, the notation $E^{w_0,n}$ is used
for a generic term of weight $w_0$, involving covariant derivatives
of $h$ of order at most $n$, satisfying certain restrictions:
$$(E^{w_0,n})\quad w[T]=w_0,\quad j_i\leq n,\quad p_i=1\mbox{ if
}j_i=n,\quad p_i=1\mbox{ or }2\mbox{ if }j_i=n-1\mbox{ and }n\geq
1.$$ The symbol $\tilde{E}^{w_0,n}$ denotes the larger space where
we drop the restrictions:
$$(\tilde{E}^{w_0,n})\quad w[T]=w_0,\quad j_i\leq n.$$
(Sometimes the same notation is used for the real vector space
spanned by terms of the given type). For example, above we checked
that:
$$(\partial_t-\Delta)|h|^2=-2|\nabla h|^2+E^{4,1},$$
$$(\partial_t-\Delta)|\nabla h|^2=-2|\nabla^2h|^2+E^{6,2}.$$
These symbols have some useful properties. For example, it is easy
to prove by induction that:
$$\nabla (E^{n+3,n+1})\subset E^{n+4,n+2},\quad n\geq 0.$$
In the proof of this the following easily verified fact (for $n>1$)
is used:
$$E^{n+3,n+1}=(\nabla^{n+1}h)*h+(\nabla^nh)*[\nabla
h+h^{(2)}]+\tilde{E}^{n+3,n-1}.$$ The property verified above for
low $n$ holds in general:\vspace{.2cm}

\textbf{Lemma 13.3.} For $n\geq 0$:
$$(\partial_t-\Delta)|\nabla^nh|^2=-2|\nabla^{n+1}h|^2+E^{2n+4,n+1}.$$

\emph{Proof} (for $n\geq 2$.) With the natural multi-index notation:
$$\partial_t|\nabla h|^2=2\langle
\partial_t(\nabla^nh),\nabla^nh\rangle+\partial_t(g^{IJ}g^{pr}g^{qs})(\nabla^n_Ih)_{pq}(\nabla^n_Jh)_{rs},
\quad |I|=|J|=n,$$ where, using the Codazzi identity and the
curvature tensor repeatedly:
$$\partial_t(g^{IJ}g^{pr}g^{qs})(\nabla^n_Ih)_{pq}(\nabla^n_Jh)_{rs}=(n+2)(\partial_tg^{ij})
\langle \nabla_i\nabla^{n-1}h,\nabla_j\nabla^{n-1}h\rangle$$
$$+(\partial_tg^{ij})Rm[\nabla^{n-2}h]_i*Rm[\nabla^{n-2}h]_j+(\partial_tg^{ij})(\nabla^nh)_i*Rm[\nabla^{n-2}h]_j.$$
Recall (appendix) $\partial_tg^{ij}=\nabla h+h^{(2)}$ and $Rm=h*h$.
Thus this term has the form:
$$(\nabla^nh)^{(2)}*(\nabla h+h^{(2)})+(\nabla^{n-2}h)^{(2)}*(\nabla
h+h^{(2)})*h^{(4)}+(\nabla^nh)*(\nabla^{n-2}h)*(\nabla
h+h^{(2)})*h^{(2)},$$ so it is in $E^{2n+4,n+1}$.

Turning to space derivatives, we have (as for $n=1$):
$$\Delta|\nabla^nh|^2=2\langle \Delta(\nabla^nh),\nabla^nh\rangle
+2|\nabla^{n+1}h|^2,$$ and therefore:
$$(\partial_t-\Delta)|\nabla^nh|^2=-2|\nabla^{n+1}h|^2+2\langle
(\partial_t-\Delta)(\nabla^nh),\nabla^nh\rangle+E^{2n+4,n+1}.$$
\vspace{.2cm}

\emph{Claim.} $(\partial_t-\Delta)[\nabla^nh]\in E^{n+3,n+1}$
($n\geq 0$).\vspace{.2cm}

The claim clearly implies the conclusion of the lemma (bearing in
mind the expression seen above for a general term in $E^{n+3,n+1}$).
Of course, we prove it by induction on $n$, the cases $n=0,1$ having
already been checked:
$$(\partial_t-\Delta)h=H\nabla_{\omega}h+T+h^{(3)}\in E^{3,1},\quad
(\partial_t-\Delta)(\nabla
h)=\nabla[(\partial_t-\Delta)h]+E^{4,2}\in E^{4,2}.$$\vspace{.2cm}

\emph{Proof of claim.} The claim would follow inductively from:
$$(\partial_t-\Delta)[\nabla^{n+1}h]=\nabla[(\partial_t-\Delta)(\nabla^nh)]+E^{n+4,n+2},$$
since $\nabla E^{n+3,n+1}\subset E^{n+4,n+2}.$ \vspace{.2cm}

For the time derivative part, we have, for any multi-index $iI$ of
length $i+1$, $|I|=n$:
$$\begin{array}{ll}\partial_t[\nabla^{n+1}h]_{iI}&=\partial_t[\partial_i(\nabla^nh[\partial_I]))
-(\nabla^nh)(\nabla_i\partial_I)]\\
&=\partial_i(\partial_t(\nabla^nh[\partial_I]))-\partial_t(\nabla^nh[\nabla_i\partial_I])\\
&=\nabla_i(\partial_t(\nabla^nh))[\partial_I]+\partial_t(\nabla^nh)(\nabla_i\partial_I)
-\partial_t(\nabla^nh)(\nabla_i\partial_I)-\nabla^nh[\partial_t(\nabla_i\partial_I)].\end{array}.$$
For a multi-index $I=i_1\ldots i_n$ of length $n$, denote by
$I^{k}_p$ the multi-index of length $n$ obtained from $I$ by setting
its $k^{th.}$ entry $i_k$ equal to $p$. It is then clear that:
$$\partial_t(\nabla_i\partial_I)=\sum_{k=1}^n\sum_p(\partial_t\Gamma_{ii_k}^p)\partial_{I^{k}_p}.$$
In symbolic notation, the preceding calculation is summarized as:
$$\partial_t[\nabla^{n+1}h]=\nabla(\partial_t\nabla^nh)+(\nabla^nh)*(\partial_t\Gamma).$$
Since $\partial_t\Gamma\in E^{3,2}$, this says:
$$\partial_t[\nabla^{n+1}h]=\nabla(\partial_t\nabla^nh)+E^{n+4,n+2}.$$

Covariant derivatives in space may be dealt with in the usual way.
Writing a multi-index $I$ of length $n+1$ as $I=iI', |I'|=n$, we
have for first-order derivatives:
$$\begin{array}{ll}\nabla_k(\nabla_I^{n+1}h)=\nabla_k(\nabla_{iI'}^{n+1}h)&=\nabla_k(\nabla_i(\nabla_{I'}^nh))
-\nabla_k(\nabla^nh(\nabla_i\partial_{I'}))\\
&=\nabla_i(\nabla_k(\nabla_{I'}^nh))-\nabla_k(\nabla^nh(\nabla_i\partial_{I'}))+Rm_{ik}[\nabla^n_{I'}h].\end{array}$$
And for second-order covariant derivatives:
$$\begin{array}{ll}\nabla_l(\nabla_k(&\nabla_I^{n+1}h))=\nabla_l(\nabla_i(\nabla_k(\nabla^n_{I'}h)))+
\nabla_l(Rm_{ik}[\nabla_{I'}^nh])-\nabla_l(\nabla_k(\nabla^nh(\nabla_i\partial_{I'})))\\
&=\nabla_i(\nabla_l(\nabla_k(\nabla^n_{I'}h)))+Rm_{il}[\nabla_k(\nabla_{I'}^nh)]
+\nabla(Rm*\nabla^nh)+\nabla^2(\nabla^nh*h)\\
&=\nabla_i(\nabla^2_{l,k}(\nabla_{I'}^nh))+\nabla_i(\nabla_{\nabla_l\partial_k}\nabla_{I'}^nh)+Rm*\nabla^{n+1}h+
\nabla(Rm*\nabla^nh)+\nabla^2(\nabla^nh*h).\end{array}$$
$$\begin{array}{ll}\nabla^2_{l,k}(\nabla_I^{n+1}h)&=\nabla_i(\nabla^2_{l,k}(\nabla^n_{I'}h))
-\nabla_{\nabla_l\partial_k}(\nabla_I^{n+1}h)+\nabla(\Gamma*\nabla^{n+1}h)\\&+Rm*\nabla^{n+1}h+
\nabla(Rm*\nabla^nh)+\nabla^2(\nabla^nh*h)\\
&=\nabla_i(\nabla^2_{l,k}(\nabla^n_{I'}h))+\Gamma*\nabla^{n+2}h+\nabla(\Gamma*\nabla^{n+1}h)\\&+Rm*\nabla^{n+1}h+
\nabla(Rm*\nabla^nh)+\nabla^2(\nabla^nh*h) .\end{array}$$ Taking
traces with $g^{kl}$ and using
$$Rm=h^{(2)},\quad \nabla Rm=\nabla h*h,\quad \Gamma=h\omega,\quad \nabla
\Gamma=\nabla h+h^{(2)},$$ it follows easily that:
$$\Delta(\nabla^{n+1}h)=\nabla(\Delta(\nabla^nh))+E^{n+4,n+2},$$
and therefore:
$$(\partial_t-\Delta)[\nabla^{n+1}h]=\nabla[(\partial_t-\Delta)(\nabla^nh)]+E^{n+4,n+2},$$
proving the claim and the lemma.\vspace{.2cm}

The analogue of Lemma 13.2 for higher covariant derivatives of $h$
follows easily from these remarks.\vspace{.2cm}

\textbf{Lemma 13.4.} For a solution $w\in
C^{n+5,[\frac{n+5}2]+1}(E)$, assume we have a uniform bound for $h$
and its first $n$ covariant derivatives: $|\nabla^jh|\leq a_j$ in
$E$, $j=0,\ldots,n$. Then we may find constants $\alpha>0,C>0$
depending only on the dimension and the $a_j$, so that the function:
$$f_{n+1}(x,t)=\alpha |\nabla^{n+1} h|^2+|\nabla^nh|^2$$
is a subsolution in $E$: $(\partial_t-\Delta)f_{n+1}\leq
C$.\vspace{.2cm}

\emph{Proof.} In the proof we denote by $C_n$ a generic positive
constant depending only on dimension and the $a_j, j=0,\ldots,n$. We
have:
$$(\partial_t-\Delta)|\nabla^nh|^2=-2|\nabla^{n+1}h|^2+E^{2n+4,n+1},\quad
(\partial_t-\Delta)|\nabla^{n+1}h|^2=-2|\nabla^{n+1}h|^2+E^{2n+6,n+2},$$
where:
$$E^{2n+4,n+1}=\nabla^{n+1}h*\nabla^nh*h+(\nabla^nh)^{(2)}*\tilde{E}^{2,1}
+(\nabla^nh)*\tilde{E}^{n+3,n-1}+\tilde{E}^{2n+4,n-1},$$
$$E^{2n+6,n+2}=\nabla^{n+2}h*\nabla^{n+1}h*h+(\nabla^{n+1}h)^{(2)}*\tilde{E}^{2,1}
+(\nabla^{n+1}h)*\tilde{E}^{n+4,n}+\tilde{E}^{2n+6,n}.$$ This
implies:
$$(\partial_t-\Delta)|\nabla^nh|^2\leq-2|\nabla^{n+1}h|^2+C_n|\nabla^{n+1}h|+C_n,$$
$$(\partial_t-\Delta)|\nabla^{n+1}h|^2\leq-2|\nabla^{n+2}h|^2+C_n|\nabla^{n+2}h||\nabla^{n+1}h|+C_n(|\nabla^{n+1}h|^2+
|\nabla^{n+1}h|+1).$$ It is easy to see from these inequalities that
we may choose $\alpha$ sufficiently small so that the conclusion of
the lemma will hold.\vspace{.4cm}

\textbf{14. H\"{o}lder gradient estimate for the second fundamental
form.}\vspace{.2cm}

\emph{Notational remark:} In this section, parabolic H\"{o}lder
spaces are denoted by a single superscript (that is, the space
previously denoted $C^{2+\alpha,(1+\alpha)/2}$ is denoted
$C^{2+\alpha}$ in this section). Capital $X,Y$ etc. denote general
points in the spacetime domain $E$. (This follows the notation used
in \cite{Lieberman}.)\vspace{.2cm}

A continuation criterion for the solution $w(y,t)$ in $E^T$ in terms
of a bound on the norm $|h|_g$ of the second fundamental form would
follow from an a-priori $C^{3+\delta}(E^T)$ bound on a solution,
assuming $|h|_g\leq a_0$ in $\bar{E}^T$; equivalently, from a global
a priori H\"{o}lder gradient bound $|\nabla h|_{\delta}\leq M$ in
$\bar{E}^T$ (for suitably controlled $M$). In this section we show
how such a bound follows from the a-priori estimates of linear
parabolic theory applied to the evolution equations for $v$, $H$,
and the Weingarten operator, under an additional
hypothesis.\vspace{.2cm}

Assuming $w\in C^{2+\delta}(E^T)$ is a solution, satisfying in
addition $|h|_g\leq a_0$ in $E^T$, we already observed the maximum
principle implies bounds:
$$0\leq w\leq w_0,\quad 1\leq v\leq\bar{v}\mbox{ in }\bar{E}^T,$$
depending only on the initial data and $\beta$ (we assume $w\geq 0$,
at $t=0$, vanishing only on $\partial D_0$.) In particular, $g$ is
uniformly equivalent to the euclidean metric on $E^T$. In this
section, bounds depending on $a_0$, $\bar{v}$ and the initial data
will be denoted generically by a constant $ M>0$ (dependence on
$\beta$ will not be recorded explicitly). The bound on $h$ implies a
uniform $C^2$ bound for the spacetime domain $E^T$, which we can
express in terms of a diffeomorphism $\Phi:D_0\times
[0,T]\rightarrow E^T$ by: $$|\Phi|_{C^2}\leq M.$$

We will also need to assume a uniform gradient bound on the boundary
for the second fundamental form:
$$|(\nabla_{\tau} h)(\tau,\tau)|\leq a_1 \quad \forall \tau\in
T\partial D(t),|\tau|=1.$$ Estimates depending $a_0$, $a_1$,
$\bar{v}$ and the initial data will be given in terms of constants
denoted generically by $M_1$.\vspace{.2cm}

In fact $E^T$ is a bounded domain in $\mathbb{R}^n\times [0,T]$ of
class $C^{2+\delta}$, with bounds controlled by $M_1$. (This
statement includes some regularity in $t$, so it is not immediate
from the uniform bound assumed for $\nabla^{tan}h^{tan}$ on
$\partial_lE$). To see this, consider the equation satisfied by
$w_k=\partial_kw$, written in `divergence form' with Dirichlet
boundary conditions:
$$\left \{
\begin{array}{l}\partial_tw_k-\partial_i(g^{ij}\partial_jw_k)=g^k:=(\partial_kg^{ij})w_{ij}-(\partial_ig^{ij})\partial_jw_k\\
{w_k}_{|\partial_lE}:=\varphi^k=(\beta_0/\beta)n^k,\quad
{w_k}_{|t=0}=\partial_kw_0\end{array}\right.$$ Assuming
$\partial_kw\in C^{1+\delta}(E)$, the following estimate holds
(\cite{Lieberman}, thm 4.27):
$$|w_k|_{1+\delta}\leq
C(\sup_E|w_k|+||g^k||_{1,n+1+\delta}+|\varphi^k|_{1+\delta;\partial_lE}+|\partial_kw_0|_{1+\delta;D_0}).$$
Here $||g^k||_{1,n+1+\delta}$ is the norm in the spacetime Morrey
space $L^{1,n+1+\delta}(E)$:
$$||g^k||_{1,n+1+\delta}=\sup_{Y\in
E,r<diam(E)}(r^{-(n+1+\delta)}\int_{E[Y,r]}|g^n|dX).$$ In the
present case this can easily be estimated, since:
$$|\partial_kg^{ij}|=|h^i_k\omega^j+h_k^j\omega^i|\leq M,\quad
|\partial_jw_k|\leq \bar{v}a_0\leq M\Rightarrow |g^k|\leq M,$$ and
$|E[Y,r]|\leq Cr^{n+2}$, while $\delta\in (0,1)$. Thus
$||g^k||_{1,n+1+\delta}\leq M$.\vspace{.2cm}

Since $|\nabla_{\tau}(\nabla_{\tau}n)|\leq
c(|(\nabla_{\tau}h)(\tau,\tau)|+|h|)\leq M_1$, it follows that $n$
is $C^2$ in space variables on $\partial_lE$. On the other hand,
$Dw=\omega/\beta$ on $\partial D(t)$, and $\omega$ is a solution of
$\partial_t\omega^k=tr_gD^2\omega^k+|h|^2_g\omega^k$, hence $n$ is
also $C^1$ in time on $\partial_lE$. We conclude
$|\varphi^k|_{1+\delta;\partial_lE}\leq(\beta_0/\beta)|n|_{1+\delta;\partial_lE}\leq
M_1$.\vspace{.2cm}

Therefore we have $|Dw|_{1+\delta}\leq M_1$, and $|w|_{2+\delta}\leq
M_1$ (note that $C$ depends on $|g^{ij}|_{C^{\delta}}$ and other
constants also controlled by $M$.) In particular, $E^T$ is a
$C^{2+\delta}$ domain with chart constants controlled by $M_1$. (In
fact, in a neighborhood of any point $P\in \partial_lE$ with
$\partial_{y_2}w\neq 0$, a boundary chart $\Psi$ is given by
$\Psi(y_1,y_2,t)=(y_1,w(y,t),t)$.)\vspace{.2cm}

The first-order term in the evolution equation for $h$ (or for the
Weingarten operator) involves $DH$; hence the next step is to obtain
a global gradient bound $|DH|_{1+\alpha}\leq M_1$ in $\bar{E}^T$.
The mean curvature satisfies the `divergence form' equation with
Neumann boundary conditions:
$$\left \{ \begin{array}{l}
\partial_tH-\partial_j(g^{ij}(X)H_{i})+\partial_j(g^{ij})(X)\partial_iH-c(X)H=0,\\
 d_nH=(\beta^2/\beta_0)Hh_{nn}:=\psi\mbox{ on }\partial D(t),\quad
H_{|t=0}=H_0,
\end{array}\right .$$
where:
$$ c:=|h|^2_g-h^2(\omega,\omega)+Hh(\omega,\omega).$$
Then, with the regularity conditions for the domain and the
coefficients:
$$\partial_lE\in C^{1+\delta},\quad n\in
C^{\delta}(\partial_lE),\quad \partial_j(g^{ij})\in
L^{1,n+1+\delta}(E),\quad c\in L^{1,n+1+\delta}(E),$$ and assuming
$H\in C^{1+\delta}(E)$, or $w\in C^{3+\delta}(E)$, we have the
bound:
$$|H|_{1+\delta;\bar{E}}\leq
C(\sup_E|H|+|\psi|_{\delta;\partial_lE}+|H_0|_{1+\delta;D_0}).$$ As
noted earlier:
$$||\partial_jg^{ij}||_{1,n+1+\delta}+||c||_{1,n+1+\delta}\leq M,$$
hence $C$ is controlled by $M$. In addition, $|w|_{2+\delta}\leq
M_1$ implies $|h|_{\delta}\leq M_1$, and hence:
$|\psi|_{\delta;\partial_lE}\leq M_1$. We conclude
$|H|_{1+\delta}\leq M_1$, and state it as a lemma.\vspace{.2cm}

\textbf{Lemma 14.1.} Let $w\in C^{3+\delta}(E^T)$ be a classical
solution of graph mean curvature motion, with contact and
constant-angle boundary conditions. Assume (i) $|h|_g\leq a_0$ on
$\partial_lE$; (ii) $|(\nabla_{\tau}h)(\tau,\tau)|\leq a_1$ on
$\partial_lE$. Then we have the global gradient bound for $H$:
$$\sup_{\bar{E}^T}|DH|_{\delta}\leq M_1,$$
for a constant $M_1$ depending on $\delta,\bar{v},a_0,a_1$ and the
initial data $w_0$.\vspace{.2cm}

\textbf{Corollary 14.2.} Under the same hypotheses as Lemma 14.1, we
have a global gradient bound:
$$\sup_{\bar{E}}|\nabla h|_g\leq M_1,$$
for a positive constant $M_1$ depending on $\delta,\bar{v},a_0,a_1$
and the initial data $w_0$.\vspace{.2cm}

\emph{Proof.} The bound on the components $(\nabla_nh)(\tau,\tau)$
and $(\nabla_nh)(n,n)$ on the lateral boundary $\partial_lE$ follows
immediately from the expressions in section 11. The bound on
$(\nabla_{\tau}h)(\tau,\tau)$ over $\partial_lE$ is hypothesized,
and then the bound on the remaining component
$(\nabla_{\tau}h)(n,n)$ follows from the global gradient bound
$|DH|\leq M$ implied by Lemma 14.1. Thus $|\nabla h|\leq M_1$ on
$\partial_lE$, and then the global bound follows from Lemma 13.2 and
the maximum principle.\vspace{.2cm}

To improve the conclusion of corollary 14.2 to  a H\"{o}lder
gradient bound, it is natural to consider the evolution equation for
$h$, with the Neumann-type boundary conditions derived in section
11. One is then faced with the problem that those boundary
conditions do not control components such as
$(\nabla_{\tau}h)(\tau,\tau)$ on $\partial_tE$. So as a preliminary
step we consider the evolution equation for $v$, which has the
advantage that the boundary values are constant. Written in linear
form, we have:
$$\left \{
\begin{array}{l}\partial_tv-g^{ij}(X)v_{ij}+b^i(X)\partial_iv+c(X)v=0\\
v_{|\partial_tE}=\frac 1{\beta},\quad v_{|t=0}=v_0,
\end{array}\right .$$
where:
$$g^{ij}(X)=\delta_{ij}-\frac{w_iw_j}{1+|Dw|^2}(X),\quad
b^i(X)=\frac{2g^{ij}w_kw_{kj}}{1+|Dw|^2}(X),\quad c(X)=|h|^2_g(X).$$
We clearly have $g^{ij}\in C^{\delta}$ (since $Dw\in C^{\delta}$),
as well as $b^i,c\in C^{\delta}$ (since $h\in C^{\delta}$) and
$\partial_lE\in C^{2+\delta}$, with bounds controlled by $M_1$ in
all cases as observed earlier. Therefore assuming $v\in
C^{2+\delta}$ (equivalently, $w\in C^{3+\delta}$) we have the bound:
$$|v|_{2+\delta;\bar{E}}\leq C(\sup_E v+\frac 1{\beta}),$$
with $C$ controlled by $M_1$. Thus $|D^2v|_{\delta}\leq M_1$.
Recalling $v^{-1}\partial_iv=h(\partial_i,\omega)$, this implies:
$$|(\nabla_{\tau}h)(n,n)|_{\delta;\partial_lE}=|(\nabla_{n}h)(\tau,n)|_{\delta;\partial_lE}\leq
M_1,\quad \forall \tau\in T\partial D(t),|\tau|=1.$$ Since
$H=\beta^2h_{nn}+ h(\tau,\tau)$ on $\partial_lE$, it follows from
Lemma 14.1 that we also have
$|(\nabla_{\tau}h)(\tau,\tau)|_{\delta;\partial_lE}\leq M_1$. For
the remaining components of $\nabla h$, this bound follows directly
from the boundary conditions:
$$|(\nabla_{n}h)(\tau,\tau)|_{\delta;\partial_lE}+|(\nabla_{n}h)(n,n)|_{\delta;\partial_lE}\leq
M_1.$$\vspace{.2cm}

 Now consider the evolution of the components of the
Weingarten operator, written in divergence form with Neumann
boundary conditions:
$$\left \{ \begin{array}{l}
\partial_th_j^k-\partial_i(g^{il}\partial_lh_j^k)=f_j^k,\mbox{ in }E^T,\quad
f_j^k:=H_jh_l^k\omega^l-H_lh_j^l\omega^k+h^{(3)k}_j-(\partial_ig^{il})(\partial_lh_j^k)\\
d_n(h_j^k)=\varphi_j^k\mbox{ on }\partial_lE,\quad
{h_j^k}_{|t=0}=h_{j0}^k
\end{array}\right .$$
The same theorem quoted above gives the estimate (assuming $h_j^k\in
C^{1+\delta}$, or $w\in C^{3+\delta}$):
$$|h_j^k|_{1+\delta;\bar{E}}\leq
C(\sup_{E}|h_j^k|+||f_j^k||_{1,n+1+\delta}+|\varphi_j^k|_{\delta;\partial_lE}+|h_{0j}^k|_{1+\delta;D_0}).$$
Note that:
$$d_n(h_j^k)=g^{ik}(\nabla_nh)_{ij}=\beta^2(\nabla_nh)(n,\partial_j)n^k+(\nabla_nh)(\tau,\partial_j)\tau^k\mbox{
on }\partial_lE.$$ From this and the above discussion it follows
that $|\varphi_j^k|_{\delta;\partial_lE}\leq M_1$. The bound
$||f_j^k||_{1,n+1+\delta}\leq M_1$ follows from Lemma 14.1 and
Corollary 14.2.  We conclude $|h_j^k|_{1+\delta;\bar{E}}\leq M_1$.
The $1+\delta$ estimate for $h_j^k$ clearly implies the following
lemma.\vspace{.2cm}

\textbf{Lemma 14.3.} Let $w\in C^{3+\delta}(E^T)$ be a classical
solution of graph mean curvature motion, with contact and
constant-angle boundary conditions. Assume (i) $|h|_g\leq a_0$ on
$\partial_lE$;(ii)$|(\nabla_{\tau}h)(\tau,\tau)|_{\partial_lE}\leq
a_1$. Then we have the global H\"{o}lder gradient bound for $h$:
$$|\nabla h|_{\delta;\bar{E}^T}\leq M_1,$$
for a constant $M_1$ depending on $\delta,\bar{v},a_0,a_1$ and the
initial data $w_0$.\vspace{.2cm}

\emph{Remark.} This is clearly equivalent to a global `a priori'
$C^{3+\delta}$ bound for $w$ on $\bar{E}^T$, $|w|_{3+\delta}\leq
M_1$.\vspace{.2cm}

Lemma 14.3 is the main step in the derivation of a `continuation
criterion' for this flow.\vspace{.2cm}

\textbf{Proposition 14.4.} Assume the maximal existence time
$T_{max}$ is finite. Then (for $n=2$, in the concave case):
$$\limsup_{t\rightarrow T_{max}}\sup_{\partial
D(t)}(|h|_g+|(\nabla_{\tau}h)(\tau,\tau)|)=\infty.$$\vspace{.2cm}

\emph{Proof.} For $w_0\in C^{3+\bar{\alpha}}(D_0)$ satisfying the
contact angle condition (with $\bar{\alpha}\in (0,1)$ arbitrary) and
$\alpha=\bar{\alpha}^2$, Theorem 8.1 yields a unique solution
$F=[u,\varphi]$ of m.c.m. with contact angle/orthogonality boundary
conditions in a maximal time interval $[0,T_{max}]$, with $F\in
C^{2+\alpha}(Q_0^{T_{max}})$, $Q_0^{T_{max}}=Q\times [0,T_{max})$;
this is also the unique solution in $F\in
C^{2+\delta^2}(Q_0^{T_{max}})$, where $\delta=\alpha^2$. Then
$w=u\circ \varphi^{-1}\in C^{2+\delta}(E^{T_{max}})$ is a solution
of graph m.c.m., which for any $t_0>0$ is in
$C^{3+\delta}(E_{t_0}^{T_{max}})$. By contradiction, assume
$|h|_g+|(\nabla_{\tau}h)(\tau,\tau)|$ is bounded in $E_{t_0}^T$ for
any $T<T_{max}$ (with bound independent of $T$). Then lemma 14.3
applies, giving an a-priori bound $|\nabla
h|_{\delta;\bar{E}_{t_0}^T}\leq M_1$, for $T$ arbitrarily close to
$T_{max}$. In particular, $|w(\cdot,T)|_{C^{3+\delta}(D(T))}\leq
M_1$, and for $T$ sufficiently close to $T_{max}$ we can use theorem
8.1 again (with initial data $w(\cdot,T)$) to find a solution
$F'=[u',\varphi']\in C^{2+\delta^2}(Q_0^{T'})$ (where $T'>T_{max}$),
extending $F$. This contradicts the maximality of
$T_{max}$.\vspace{.5cm}

\textbf{15. Behavior at the extinction time.}\vspace{.2cm}

In this section we consider the behavior of $\Sigma_t$ as $t$
approaches the maximal existence time $T$, in the concave case. We
also assume $H\leq H_0<0$ at $t=0$, so $T$ is finite. Let
$K_t\subset {\mathbb R}^{n+1}$ be the compact convex set bounded by
$\Sigma_t$. Since $H\leq 0$, $\{K_t\}$ is a decreasing family, and
the intersection:
$$K_T=\bigcap_{0\leq t<T}K_t\subset {\mathbb R}^{n+1}$$ is compact,
convex and non-empty. It turns out $K_T$ has zero $(n+1)$-volume. In
this section we use the support function to show this when $n=2$
(following the argument in \cite{Stahl}), under the assumption that
there is no `gradient breakdown'.

Assume the origin $0\in {\mathbb R}^{n+1}$ is a point of $K_T$. The
support function of $K_t$ (with respect to this origin) is the
function $p(\cdot,t)$ on $D(t)$:
$$p(y,t)=\langle G(y,t),N(y,t)\rangle, \quad G(y,t)=[y,w(y,t)].$$
Since $K_t$ is convex, $p>0$ in $D(t)$; the evolution equations and
boundary conditions for $p$ are easily computed: from $L[G]=0$ and
$L[N]=|h|_g^2N$, we have:
$$L[p]=\langle L[G],N\rangle+\langle G,L[N]\rangle-2g^{kl}\langle
\partial_kG,\partial_lN\rangle$$
$$=|h|^2_gp+2H,$$
and since $\langle d_nG,N\rangle=0$: $${p_n}_{|\partial
D(t)}=\langle G,d_nN\rangle=-A(G^T,N),$$ where, with $y^T:=y-(y\cdot
n)n\in T_y\partial D(t)$, the tangential component $G^T:=G-\langle
G,N\rangle N$ is easily seen to be, at $\partial D(t)$:
$$G^T=\frac 1{v^2}[w_n^2y^T+y,0].$$
Since $h(y^T,n)=0$ at $\partial D(t)$, this implies
$A(G^T,n)=\beta^2(y\cdot n)h(n,n). $ Note that $p(y)=-\beta_0
(y\cdot n)$ on $\partial D(t)$, so we have:
$${p_n}_{|\partial D(t)}=\frac{\beta^2}{\beta_0}ph_{nn},$$
reminiscent of the boundary condition for $H$. Note also that we
have the upper bound:
$$p\leq||G||\leq \max_{D(0)}||G_0||:=p_0,$$
since the $K_t$ are `nested'.
 \vspace{.2cm}

\textbf{Proposition 15.1} ($n=2$)  Assume we know that
$$\limsup_{t\rightarrow T}(\sup_{y\in \partial D(t)}|h|_g)=\infty$$
at the maximal existence time $T$. Then: $$\liminf_{t\rightarrow
T}(\inf_{y\in D(t)}p(y,t))=0.$$\vspace{.2cm}

\emph{Proof.} Reasoning by contradiction, assume $p>2\delta >0$ for
$t\in [0,T).$ We claim that this implies an upper bound for $|H|$
(and hence for $|h|$, since $|h|^2\leq nH^2$), contradicting the
fact that $\limsup_{t\rightarrow T}\sup_{\Gamma_t}|h|=\infty$.
\vspace{.2cm}

To prove the claim, consider the function:
$$f(y,t):=\frac{|H|}{p-\delta}=-\frac{H}{p-\delta}.$$
Using the evolution equations and boundary conditions for $H$ and
$p$ we find (with $\hat{\omega}:=\omega/|\omega|_g$, see Remark 15.1
below):
$$L[f]=f(-\delta
|h|^2_g+2pf+|\omega|^2_g(h^2(\hat{\omega},\hat{\omega})-Hh(\hat{\omega},\hat{\omega}))
-\frac 2{p-\delta}g^{kl}\partial_kf\partial_lp,$$
$${f_n}_{|\partial D(t)}=-\delta
\frac{\beta^2}{\beta_0}h_{nn}\frac{|h|^2}{(p-\delta)^2}\geq 0.$$
Since $h|_g^2\geq (1/n)H^2=(1/n)f^2(p-\delta)^2$:
$$L[f]\leq -\frac{(p-\delta)^2\delta}nf^2+2pf+
|\omega|^2_g(h^2(\hat{\omega},\hat{\omega})-Hh(\hat{\omega},\hat{\omega}))
-\frac 2{p-\delta}\langle \nabla f,\nabla p\rangle_g.$$ Now recall
(Remark 12.1) that if $n=2$:
$$h^2(\hat{\omega},\hat{\omega})-Hh(\hat{\omega},\hat{\omega})=-\Delta\leq
0,$$ so:
$$L[f]\leq f(-\frac{(p-\delta)^2\delta}nf^2+2pf)-\frac
2{p-\delta}\langle \nabla f,\nabla p\rangle_g.$$ Assume $\delta>0$
is so small that $\sup_{D(0)}f_{|t=0}<\frac{2np_0}{\delta^3}.$ We
claim this persists for all $t\in [0,T)$. If not, assume
$f(y_0,t_0)=2np_0/\delta^3$, with $t_0>0$ smallest possible and
$y_0$ a local max of $f(\cdot, t_0)$. Since $f_n\geq 0$ at $\partial
D(t_0)$, $z_0=(y_0,t_0)$ can't be a boundary point of $E$ (by the
boundary point lemma). Since $y_0\in \partial D(t_0)$ is an interior
point, we have $L[f]_{|z_0}\geq 0$ and $\nabla f(z_0)=0$, hence:
$$\frac{(p-\delta)^2\delta}nf(z_0)\leq 2p(z_0),\mbox{ or
}f(z_0)\leq \frac{2np(z_0)}{\delta(p-\delta)^2}\leq
\frac{2np_0}{\delta(p-\delta)^2},$$ which is not possible since
$p-\delta >\delta$. Thus $f(y,t)<4p_0/\delta^3$ in $E$, which
implies the bound:
$$|H|\leq \frac {4p_0^2}{\delta^3}$$ for $t\in [0,T)$,
contradicting the maximality of $T$.\vspace{.3cm}

\emph{Remark 15.1.} Consider the vector fields in $D(t)\subset
\mathbb{R}^2$:
$$\omega=\frac 1v[w_1,w_2], \quad
\tilde{\omega}=v\omega^{\perp}=[-w_2,w_1].$$ It is easy to verify
the following:
$$\langle \omega,\tilde{\omega}\rangle_g=0,\quad
|\omega|^2_g=|\tilde{\omega}|^2_g=|Dw|^2_e:=w_1^2+w_2^2.$$ Thus we
may think of $\{\omega, \tilde{\omega}\}$ as a `conformal
pseudo-frame' ($\omega$ and $\tilde{\omega}$ vanish when $Dw=0$),
defined on all of $D(t)$. Moreover, at the boundary $\partial D(t)$:
$$\omega=\beta_0n,\quad
\tilde{\omega}=\frac{\beta_0}{\beta}n^{\perp}:=\frac{\beta_0}{\beta}\tau,$$
where $\{\tau,n\}$ is an euclidean-orthonormal frame along
$\Gamma_t$. Thus $\omega$ and $\tilde{\omega}$ supply `canonical'
extensions of $n,\tau$ to the interior of $D(t)$, as uniformly
bounded vector fields.\vspace{.3cm}

\emph{Remark 15.2.} It follows from the proposition that $K_T$
cannot contain a half-ball of positive radius centered at a point of
${\mathbb R}^2$; in particular, $vol_3(K_T)=0$. Based on the
experience with curve networks \cite{LensSeminar}, one is led to
expect that $K_T$ is a point (that is, that $diam(K_T)=0$), at least
under the same assumption as the proposition (`no gradient
breakdown'). We have not been able to show this yet; existence of
self-similar solutions and comparison arguments appropriate to the
free-boundary setting appear to be needed for the usual approach to
work. \vspace{.4cm}

\textbf{16. Final comments.}\vspace{.2cm}

1. We state here the local existence theorem for configurations of
graphs over domains with moving boundaries. In this setting, a
\emph{triple junction configuration} consists of three embedded
hypersurfaces $\Sigma^1,\Sigma^2,\Sigma^3$ in $\mathbb{R}^{n+1}$,
graphs of functions $w^I$ defined over time-dependent domains
$D^1(t), D^2(t)\subset \mathbb{R}^n$ ($D^1$ covered by one graph,
$D^2$ by two graphs), satisfying the following conditions: (1) The
$\Sigma^I$ intersect along an $(n-1)$-dimensional graph $\Lambda(t)$
(the `junction'), along which the upward unit normals satisfy the
relation: $N_1+N_2=N_3$. (2) If a fixed support hypersurface
$S\subset {\mathbb R}^{n+1}$ is given (also a graph, not necessarily
connected), the $\Sigma^I$ intersect $S$ orthogonally.

Topologically, in the case of bounded domains one has the following
examples: (i) (`lens' type) 2 disks (or two annuli) covering
$D^2(t)$ and one annulus covering $D^1(t)$; (ii) (`exterior' type)
two annuli covering $D^2(t)$ and one disk covering $D^1(t)$. The
boundary component of the annuli disjoint from the junction
intersects the support hypersurface $S$ orthogonally for each
$t$.\vspace{.2cm}

 Let $\Sigma_0^I$ ($I=1,2,3$) be graphs of $C^{3+\alpha}$
functions over $C^{3+\alpha}$ domains $D_0^1,D_0^2\subset
\mathbb{R}^n$, defining a triple junction configuration and
satisfying the compatibility condition for the mean curvatures on
the common boundary $\Gamma_0$ of $D^1_0$ and $D^2_0$:
$$H^1+H^2=H^3.$$
Then there exists $T>0$ depending only on the initial data, and
functions $w^I\in C^{2+\alpha,1+\alpha/2}(Q^I)$, $Q^I\subset
{\mathbb R}^n\times [0,T)$, so that the graphs of
$w^I(.,t):D^I(t)\rightarrow {\mathbb R}$ define a triple junction
configuration for each $t\in [0,T)$, moving by mean curvature.

The proof will be given elsewhere. \vspace{.2cm}

2. An interesting issue we have not addressed here is whether one
has breakdown of uniqueness for initial data of lower regularity, or
if the `orthogonality condition' at the junction is removed. For
curve networks, non-uniqueness has been considered in
\cite{MazzeoSaez}; but neither a drop in regularity (from initial
data to solution, in H\"{o}lder spaces) nor the orthogonality
condition play a role in the case of curves.\vspace{.5cm}

\emph{Appendix 1}: \textbf{Proof of lemma 4.1.}\vspace{.2cm}

Throughout the proof, $n$ denotes the inner unit normal at $\partial
D$, extended to a tubular neighborhood $\cal{N}$ so that $D_nn=0$.
Since $D$ is  uniformly $C^{3+\alpha}$, if follows that $n\in
C^{2+\alpha}(\partial D)$, with uniform bounds. Denote by $\rho$ the
distance to the boundary (so $D\rho=n$ in $\cal{N}$). Let $\zeta\in
C^3(\bar{D})$ be a cutoff function, with $\zeta\equiv 1$ in ${\cal
N}_1\subset {\cal N}$, $\zeta \equiv 0$ in $D\setminus {\cal
N}$.\vspace{.2cm}

We find $\varphi$ of the form:
$$\varphi(x)=x+\zeta(x)f(x)n(x)$$
with $f\in C^{2+\alpha}({\cal N})$. The 1-jet conditions on
$\varphi$ at $\partial D$ translate to the conditions on $f$:
$$f_{|\partial D}=0,\quad Df_{|\partial D}=0,\quad
D^2f(n,n)_{|\partial D}={\Delta f}_{|\partial D}=h.$$ Now
use:\vspace{.2cm}

\emph{Lemma A.1.} Let $D$ be a uniformly $C^{3+\alpha}$ domain with
boundary distance function $\rho>0$. Let $h\in C^{\alpha}(\partial
D)$ be a bounded function. Then there exists an extension $g\in
C^{\infty}(D)\cap C(\bar{D})$ so that $g_{|\partial D}=h$,
$\sup_{\bar{D}}|g|\leq \sup_{\partial D}|h|$ and $\rho^2 g\in
C^{2+\alpha}(\bar{D})$.

Given this lemma, all we have to do is set $f=(1/2)\rho^2g$, which
clearly satisfies all the requirements (in particular, $\Delta f=h$
at $\partial D$.)\vspace{.2cm}

To verify that $\varphi$ is a diffeomorphism, it suffices to check
that $|\zeta fn|_{C^1}$ (in ${\cal N}\subset \{ \rho<\rho_0\}$) is
small if $\rho_0$ is small. This is easily seen:

$$|\zeta f n|_{C^0}\leq (1/2)\rho_0^2|g|_{C^0};$$
$$|D\zeta|\leq c\rho_0^{-1}\Rightarrow |fD\zeta|\leq c\rho_0
|g|_{C^0}.$$
$$|Df|\leq (1/2)\rho_0^{\alpha}||g||_{C^{2+\alpha}(\bar{D})}$$
on $\cal{N}$, since $Df\in C^{1+\alpha}(\bar{D})$ and $Df_{|\partial
D}=0$. And finally, with ${\cal A}$ the second fundamental form of
$\partial D$:
$$|Dn|\leq |{\cal A}|_{C^0}\Rightarrow
|fDn|\leq(1/2)\rho_0^2|g|_{C^0}|{\cal A}|_{C^0}.$$\vspace{.2cm}

A word about Lemma A.1. (This is probably in the literature, but I
don't know a reference.) If $D$ is the upper half-space, we solve
$\Delta g=0$ in $D$ with boundary values $h$. Then the estimate
$$[D^2(\rho^2 P*h)]^{(\alpha)}(\bar{D})\leq c
|h|_{C^{\alpha}(\partial D)}$$ follows by direct computation with
the Poisson kernel $P$; for the rest of the norm, use interpolation.
Then transfer the estimate to a general domain using `adapted local
charts', in which $\rho$ in $D$ corresponds to the vertical
coordinate in the upper half-space. (It is easy to see that at each
boundary point there is a $C^{2+\alpha}$ adapted chart, with uniform
bounds.)\vspace{.5cm}

\noindent \emph{\textbf{Appendix 2: Evolution equations for the
second fundamental form.}}\vspace{.2cm}

We consider mean curvature motion of graphs:
$$G(y,t)=[y,w(y,t)],\quad y\in D(t)\subset \mathbb{R}^n,$$
$$w_t=g^{ij}w_{ij}=vH,\quad v=\sqrt{1+|Dw|^2}.$$
In this appendix we include evolution equations for geometric
quantities, in terms of the operators:
$$\partial_t-\Delta_g,\quad \quad L=\partial_t-tr_gd^2.$$
It is often convenient to use the vector field in $D(t)$:
$$\omega:=\frac 1v Dw.$$
Since $-\omega$ is the ${\mathbb R^n}$ component of the unit normal
$N$ and $L[N]=|h|^2_gN$, we have:
$$L[\omega^i]=|h|^2_g\omega^i, \quad |h|^2_g:=g^{ik}g^{jl}h_{ij}h_{kl}.$$
Here  $h=(h_{ij})$ is the pullback to $D(t)$ of the second
fundamental form $A$:
$$h(\partial_i,\partial_j)=h_{ij}=A(G_i,G_j)=\frac
1vw_{ij}.$$

First, denoting by $\nabla$ the pullback to $D(t)$ of the induced
connection $\nabla^{\Sigma}$ (that is,
$G_*(\nabla_XY)=\nabla_{G_*X}^{\Sigma}G_*Y$ for any vector fields
$X,Y$ in $D(t)$), and using the definition:
$$\nabla^{\Sigma}_{G_i}G_j=G_{ij}-\langle G_{ij},N\rangle
N=[0,w_{ij}]-\frac
1{v^2}w_{ij}[-Dw,1]=\frac{w_{ij}}{v^2}[Dw,|Dw|^2]=\frac{w_{ij}}{v^2}G_*Dw,$$
we conclude:
$$\nabla_{\partial_i}\partial_j=\frac 1vh_{ij}Dw=h_{ij}\omega.$$

From this one derives easily a useful expression relating the
Laplace-Beltrami operator and the operator $tr_gd^2$ acting on
functions:
$$\Delta_gf=tr_gd^2f-\frac Hvw_mf_m=tr_gd^2f-Hd_{\omega}f.$$
\vspace{.2cm}

We also have, for the covariant derivatives of $h$ with respect to
the euclidean connection and to $\nabla=\nabla^g$:
$$\partial_m(h_{ij})=\nabla_mh_{ij}+[h_{jm}h_{ik}+h_{im}h_{jk}]\omega^k.$$
(Here $\nabla h$ is the symmetric $(3,0)$-tensor with components:
$\nabla_mh_{ij}=(\nabla_{\partial_m}h)(\partial_i,\partial_j)$.)\vspace{.2cm}

Iterating this and taking $g$-traces yields (using the Codazzi
identity and the easily verified relation
$\partial_i\omega^k=h^k_i:=g^{jk}h_{ij}$):
$$tr_gd^2(h_{ij})=g^{mk}\partial_m(\partial_k(h_{ij}))=g^{mk}(\nabla^2_{\partial_m,\partial_k}h)(\partial_i,\partial_j)$$
$$+H\nabla_{\omega}h_{ij}+2[h_i^k\nabla_kh_{jp}+h_j^k\nabla_kh_{ip}]\omega^p+[H_ih_{jp}+H_jh_{ip}]\omega^p$$
$$+2[h_{ip}(h^2)_{jq}+(h^2)_{ip}h_{jq}+Hh_{ip}h_{jq}]\omega^p\omega^q+2(h^3)_{ij}+2(h^2)_{ij}h(\omega,\omega).$$

Here the powers $h^2$ and $h^3$ of $h$ are the symmetric 2-tensors
defined used the metric:
$$(h^2)_{ij}:=g^{kp}h_{ik}h_{pj}=h^k_ih_{pj},\quad (h^3)_{ij}:=g^{kp}g^{lq}h_{ik}h_{pl}h_{qj}.$$
Note also that:
$$[h_i^k\nabla_kh_{jp}+h_j^k\nabla_kh_{ip}]\omega^p=\nabla_{\omega}(h^2)_{ij},$$
using the Codazzi identity.\vspace{.2cm}

\emph{Evolution equations for $h$.}

Starting from $G_t=vHe_{n+1}=H(N+\frac 1v[Dw,|Dw|^2])=HN+HG_*\omega$
and $N_t=-\nabla^{\Sigma} H-Hv^{-1}\nabla^{\Sigma}v$ (where
$\nabla^{\Sigma}f=g^{ij}f_jG_i$ and $\nabla f=g^{ij}f_j\partial_i$)
we have:
$$\partial_t(h_{ij})=\langle (HN)_{ij},N\rangle-\langle
G_{ij},\nabla^{\Sigma}H\rangle-\frac Hv\langle G_{ij},\nabla^\Sigma
v\rangle+\langle (HG_*\omega)_{ij},N\rangle.$$ Using the easily
derived facts:
$$\langle N_{ij},N\rangle=-h^2(\partial_i,\partial_j),$$
$$H_{ij}-\langle G_{ij},\nabla^{\Sigma}H\rangle=(\nabla
dH)(\partial_i ,\partial_j),$$
$$\frac 1v\langle G_{ij},\nabla^{\Sigma}
v\rangle=h(\omega,\omega)h_{ij},$$ we obtain:
$$\partial_t(h_{ij})=(\nabla
dH)(\partial_i,\partial_j)-Hh^2(\partial_i,\partial_j)-Hh(\omega,\omega)h_{ij}+
\langle (HG_*\omega)_{ij},N\rangle,$$ where:
$$\langle (HG_*\omega)_{ij},N\rangle=H_i\langle
(G_*\omega)_j,N\rangle+H_j\langle (G_*\omega)_i,N\rangle+H\langle
(G_*\omega)_{ij},N\rangle.$$
 To identify the terms, computation shows that:
$$\langle (G_*\omega)_i,N\rangle=h(\omega,\partial_i),$$
and hence, using also:
$$\nabla^{\Sigma}_{G_i}(G_*\omega)=G_*(\nabla_{\partial_i}\omega),\quad
\nabla_{\partial_i}\omega=(h_i^p+\omega^qh_{iq}\omega^p)\partial_p=\sum_ph_{ip}\partial_p,$$
we obtain (using
$\omega^k\partial_j(h_{ik})=\nabla_{\omega}h_{ij}+2h(\partial_i,\omega)h(\partial_j,\omega)$):
$$\langle
(G_*\omega)_{ij},N\rangle=\partial_j(\omega^kh_{ik})-\langle
\nabla_{G_i}^{\Sigma}(G_*\omega),\partial_jN\rangle
=h_j^kh_{ik}+\omega^k\partial_j(h_{ik})+h(\partial_j,\nabla_{\partial_
i}\omega)$$
$$=(\nabla_{\omega}h)_{ij}+(h^2)_{ij}+2h(\omega,\partial_i)h(\omega,\partial_j)+\sum_ph_{ip}h_{jp}$$
$$=(\nabla_{\omega}h)_{ij}+2(h^2)_{ij}+3h(\omega,\partial_i)h(\omega,\partial_j),$$
since
$\sum_ph_{ip}h_{jp}=(h^2)_{ij}+h(\omega,\partial_i)h(\omega,\partial_j).$
 Combining all the terms yields the result:
$$\partial_t(h_{ij})=(\nabla
dH)(\partial_i,\partial_j)+H\nabla_{\omega}h_{ij}+H_ih(\omega,\partial_j)+H_jh(\omega,\partial_i)$$
$$+H(h^2)_{ij}+3Hh(\omega,\partial_i)h(\omega,\partial_j)-Hh(\omega,\omega)h_{ij}.$$\vspace{.2cm}

From this expression and Simons' identity (in tensorial form):
$$\nabla dH=\Delta_gh+|h|^2_gh-Hh^2,$$
we obtain easily a tensorial `heat equation' for $h$:
$$[(\partial_t-\Delta_g)h]_{ij}=H\nabla_{\omega}h_{ij}+H_ih(\omega,\partial_j)+H_jh(\omega,\partial_i)+C_{ij},$$
$$C_{ij}:=|h|^2_gh_{ij}+3Hh(\partial_i,\omega)h(\partial_j,\omega)-Hh(\omega,\omega)h_{ij}.$$\vspace{.2cm}

Using the earlier computation relating $\Delta_gh$ (the tensorial
Laplacian  of $h$) and $tr_gd^2h$, we obtain from this the evolution
equation in terms of $L$:
$$L[h_{ij}]=-2[h_i^k\nabla_{\omega}h_{jk}+h_j^k\nabla_{\omega}h_{ik}]+\tilde{C}_{ij},$$
$$\tilde{C}_{ij}:=C_{ij}-2[h(\partial_i,\omega)h^2(\partial_j,\omega)+h^2(\partial_i,\omega)h(\partial_j,\omega)]
-2(h^3)_{ij}-2(h^2)_{ij}h(\omega,\omega)-2Hh(\partial_i,\omega)h(\partial_j,\omega)$$

We may also write this purely in terms of the euclidean connection
$d$:
$$L[h_{ij}]=-2[h_i^kd_{\omega}h_{jk}+h_j^kd_{\omega}h_{ik}]+\bar{C}_{ij},$$
$$\bar{C}_{ij}=C_{ij}+2[h(\partial_i,\omega)h^2(\partial_j,\omega)+h^2(\partial_i,\omega)h(\partial_j,\omega)]
-2(h^3)_{ij}-2(h^2)_{ij}h(\omega,\omega)-2Hh(\partial_i,\omega)h(\partial_j,\omega)$$\vspace{.2cm}

\emph{Time derivatives and evolution equations for $\omega$ and
$g$.}

It is sometimes convenient to use the `Weingarten operator':
$$S(X):=S(X^i\partial_i)=h^i_jX^j\partial_i.$$

The time derivative of $\omega$ is simply minus the time derivative
of the $\mathbb{R}^n$ component of $N$. In addition, one computes
easily that $\frac{\nabla v}v=S(\omega)$, so we have:
$$\partial_t\omega=\nabla H+\frac Hv \nabla v=\nabla H+HS(\omega).$$

For the metric and `inverse metric' tensors we have: from
$\partial_tg_{ij}=(w_iw_j)_t$ and $w_{it}=(vH)_i$:
$$\partial_tg_{ij}=v^2(H_i\omega^j+H_j\omega^i)+v^2H(h(\omega,\partial_i)\omega^j+h(\omega,\partial_j)\omega^i),$$
and then, using $\partial_tg^{ij}=-g^{ik}\partial_tg_{kl}g^{lj}$:
$$\partial_tg^{ij}=- [(\nabla H)^i\omega^j+(\nabla
H^j)\omega^i]-H[S(\omega)^i\omega^j+S(\omega)^j\omega^i].$$

Since we know the evolution equation of $\omega$, it is easy to
obtain that of $g^{ij}$:
$$L[g^{ij}]=-L[\omega^i\omega^j]=-L[\omega^i]\omega^j+2g^{kl}(\partial_k\omega^i)(\partial_l\omega^j)-\omega^iL[\omega^j].$$
Using $\partial_k\omega^i=h_k^i$, we find:
$$L[g^{ij}]=-2|h|^2_g\omega^i\omega^j+2(h^2)^{ij}.$$
It is also easy to see that
$\partial_kg^{ij}=-(h_k^i\omega^j+h_k^j\omega^i)$.\vspace{.2cm}

\emph{Evolution of mean curvature.}

To compute the evolution equation for $H=g^{ij}h_{ij}$, we just need
to remember $g^{ij}$ is time-dependent:
$$(\partial_t-\Delta_g)H=(\partial_tg^{ij})(h_{ij})+tr_g[(\partial_t-\Delta_g)h]
=-2h(\nabla
H,\omega)-2Hh^2(\omega,\omega)+tr_g[(\partial_t-\Delta_g)h].$$ The
result is:
$$(\partial_t-\Delta_g)H=Hd_{\omega}H+|h|^2_gH+Hh^2(\omega,\omega)-H^2h(\omega,\omega).$$
Since $L[f]=(\partial_t-\Delta_g)f-Hd_{\omega}f$ (for any $f$), we
see that the equation in terms of $L$ has no first-order terms:
$$L[H]=|h|^2_gH+Hh^2(\omega,\omega)-H^2h(\omega,\omega)$$

\emph{Remark.} One can also find $L[H]$ starting from the
expression:
$$L[g^{ij}h_{ij}]=L[g^{ij}]h_{ij}+g^{ij}L[h_{ij}]-2g^{kl}(\partial_kg^{ij})(\partial_lh_{ij}).$$
This may be used to check  the calculation.\vspace{.3cm}

\emph{Evolution of the Weingarten operator.}

The tensorial Laplacian of $S$ is the $(1,1)$ tensor $\Delta_g S$
with components $\Delta_gh^k_j$. We have:
$$\Delta_gh^k_j=g^{ik}\Delta_gh_{ij},\quad \mbox{ or }\langle
(\Delta_gS)X,Y\rangle_g=(\Delta_gh)(X,Y).$$ The evolution equation
is easily obtained:
$$(\partial_t-\Delta_g)h_j^k=(\partial_t
g^{ik})h_{ij}+g^{ik}(\partial_t-\Delta_g)h_{ij}$$
$$=H\nabla_{\omega}h_j^k+H_jh_l^k\omega^l-H_lh_j^l\omega^k+|h|^2_gh_j^k
+2HS(\omega)^kh(\omega,\partial_j)-Hh(\omega,\omega)h_j^k-Hh(S(\omega),\partial_j)\omega^k.$$
\emph{Remark:} Since the components of $\nabla S$ are given by:
$$(\nabla_{\omega}S)(\partial_j)=(\nabla_{\omega}h_j^k)\partial_k,\quad
\nabla_{\omega}h_j^k=d_{\omega}(h_j^k)+h^2(\omega,\partial_j)\omega^k-h(\omega,\partial_j)S(\omega)^k,$$
we see that upon setting $j=k$ and adding over $k$ we recover the
evolution equation for $H$.\vspace{.2cm}

The evolution equation for $h^k_j$ in terms of $L$ follows from the
calculation:
$$L[h^k_j]=L[g^{ik}]h_{ij}+g^{ik}L[h_{ij}]-2g^{mn}(\partial_mg^{ik})(\partial_nh_{ij})$$
$$=-2(\nabla_{\omega}h_m^k)h_j^m+(\partial_j|h|^2_g)\omega^k$$
$$+|h|_g^2h_j^k-Hh(\omega,\omega)h_j^k+HS(\omega)^kh(\partial_j,\omega)
+2h^3(\partial_j,\omega)\omega^k-2(h^2)^k_p\omega^ph(\partial_j,\omega).$$
Setting $j=k$ and adding over $k$, we recover the earlier expression
for $L[H]$. \vspace{.2cm}

\emph{Evolution of $|h|^2_g$.}

The fact that $g^{ij}$ is time-dependent introduces an additional
term in the usual expression:
$$(\partial_t-\Delta_g)|h|^2_g=-2|\nabla h|_g^2+2\langle
h,(\partial_t-\Delta_g)h\rangle_g+2(\partial_tg^{ij})(h^2)_{ij}.$$
Using the expressions given earlier, one easily finds:
$$(\partial_t-\Delta_g)|h|^2_g=-2|\nabla
h|_g^2+Hd_{\omega}|h|^2_g+2|h|^4_g-4Hh^3(\omega,\omega)-2H|h|^2_gh(\omega,\omega),$$
$$L[|h|^2_g]=-2|\nabla h|^2_g+2|h|^4_g-4Hh^3(\omega,\omega)-2H|h|^2_gh(\omega,\omega).$$

\end{document}